\documentclass[preprint,3p]{elsarticle}

\RequirePackage[OT1]{fontenc}
\RequirePackage{amsthm,amsmath,amssymb}
\RequirePackage{natbib}
\RequirePackage[colorlinks,citecolor=blue,urlcolor=blue]{hyperref}
\usepackage{stmaryrd}
\usepackage[sans]{dsfont}

\usepackage[english]{babel}
\usepackage{latexsym}
\usepackage{bbm}
\usepackage[mathcal]{euscript}
\usepackage{graphicx}
\usepackage{pstricks}
\usepackage{pstricks-add}
\usepackage{pst-plot}
\usepackage{enumerate}
\usepackage{multicol}
\usepackage{subfigure}
\newcommand{\e}{\varepsilon}   
\newcommand{\ow}{\overline\omega}
\newcommand{\RR}{\mathbb{R}}    
\newcommand{\CC}{\mathbb{C}}    
\newcommand{\ind}{\mathbbm{1}} 
\newcommand{\A}{\mathcal{A}}    
\newcommand{\C}{\mathcal{C}}    

\newcommand{\E}{\mathcal{E}}    
\newcommand{\F}{\mathcal{F}}    
\newcommand{\G}{\mathcal{G}} 
\newcommand{\h}{\mathbf{H}}
\renewcommand{\L}{\mathcal{L}}
\newcommand{\NN}{\mathbb{N}}

\renewcommand{\O}{\Omega}    %
\renewcommand{\o}{\omega}    %
\renewcommand{\P}{\mathcal{P}}

\newcommand{\PP}{\mathbb{P}}    
\newcommand{\EE}{\mathbb{E}}     

\newcommand{\df}{\stackrel{\mathrm{def}}=}

\newcommand{\sub}{\subseteq}    

\providecommand{\rpar}[1]{\left( #1 \right)}               
\providecommand{\kpar}[1]{\left\{ #1 \right\}}              
\providecommand{\spar}[1]{\left[  #1 \right]}              
\providecommand{\ppar}[1]{\left\langle #1 \right\rangle}     
\providecommand{\abs}[1]{\left\vert #1 \right\vert}           
\providecommand{\norm}[1]{\left\Vert #1 \right\Vert}          

\providecommand{\pd}[2] {\frac{\partial #1} {\partial #2}}  

\numberwithin{equation}{section}
\theoremstyle{plain}
\newtheorem{theorem}{Theorem}[section]
\newtheorem{lemma}[theorem]{Lemma}
\newtheorem{proposition}[theorem]{Proposition}
\newtheorem{corollary}[theorem]{Corollary}
\newtheorem{remark}{Remark}[section]
\theoremstyle{definition}

\newcommand{\vga}{\vec g}
\newcommand{\hull}{H_{\vec{g}} }
\newcommand{\sbm}{{SBM}}

\begin{document}

\begin{frontmatter}

\title{Spinning Brownian motion}

\author{Mauricio A. Duarte} \ead{mauricio.duarte@unab.cl}
\address{ Universidad Andres Bello \\ Departamento de Matem\'atica\\ Rep\'ublica 220, Piso 2\\ Santiago, Chile}

\begin{abstract}
%

We prove strong existence and uniqueness for a reflection process $X$ in a smooth, bounded domain $D$ that behaves like obliquely-reflected-Brownian-motion, except that the direction of reflection depends on a (spin) parameter $S$, which only changes when $X$ is on the boundary of $D$ according to a physical rule. The process $(X,S)$ is a degenerate diffusion.

We show uniqueness of the stationary distribution by using techniques based on excursions of $X$ from $\partial D$, and an associated exit system. We also show that the process admits a submartingale formulation and use related results to show examples of  the stationary distribution.
\end{abstract}

\end{frontmatter}

\section{Introduction}
\label{intro}

Let $D\sub\RR^n$ be a bounded $C^2$ domain, and let $B_t$ be an $n$-dimensional Brownian motion. A pair $(X_t,S_t)$ with values in $\overline D\times\RR^p$ is called \textbf{spinning Brownian motion (\sbm)} if it solves the following stochastic differential equation
\begin{align}
\label{eq:sbm}
\left\{
\begin{array}{rl}
dX_t &= \sigma(X_t)dB_t + \vec{\gamma}(X_t,S_t)dL_t, \\
dS_t &= \spar{\vec{g}(X_t)-\alpha(X_t)S_t}dL_t,
\end{array}
\right.
\end{align}
where $L_t$ is the local time for $X_t$ on the boundary of $\partial D$, and $\vec{\gamma}$ points uniformly into $D$. Our assumptions on the coefficients are as follows:
\begin{itemize}
 \item $\sigma(\cdot)$ is an $(n\times n)$-matrix valued, Lipschitz continuous function, and is uniformly elliptic, that is, there is a constant $c_1>0$ such that $\xi^T\sigma(x)\xi\geq c_1\abs{\xi}^2$ for all $\xi\in\RR^n$, and all $x\in\overline D$.
 \item $\vec\gamma(x,s) = \vec{n}(x) + \vec{\tau}(x,s)$ is defined for $x\in\partial D$ and $ s\in\RR^p$, where $\vec{n}$ is the interior normal to $\partial D$, and $\vec\tau$ is a Lipschitz vector field on $\partial D\times\RR^p$ such that $\vec{n}(x)\cdot\vec\tau(x,s)=0$ for all $x\in\partial D$ and $s\in\RR^p$,
 \item $\vec{g}(\cdot)$ is a Lipschitz vector field on $\partial D$ with values in $\RR^p$. The function $\alpha(\cdot)$ is defined on $\partial D$, it is assumed to be Lipschitz, with values in a compact interval $[ \alpha_0,\alpha_1]$ for some $0<\alpha_0<\alpha_1<\infty$.
\end{itemize}

The process $X_t$ behaves just like a Brownian diffusion inside $D$, and is reflected instantaneously in the direction $\vec\gamma=\vec{n}+\vec\tau$ when it hits the boundary. The challenge is that the direction of reflection depends on the multidimensional parameter $S_t$, which is updated every time the main process $X_t$ hits the boundary of $D$.

This type of process arises naturally from a physical model that might be useful for applications: consider a small ball that spins and moves around a planar box following a Brownian path. On the boundary of the box we put tiny wheels which rotate at different speeds, modifying the spin of the ball as well as pushing it in a certain (non-tangential) direction. In this context, it is natural to think of the boundary wheels as an external forcing system that is not affected by the hitting of the ball: every wheel on the boundary rotates at a speed dependent only on its position. The position $X_t$ of the particle at time $t$ is described by the first equation in \eqref{eq:sbm}, in which the direction of the boundary push $\vec\gamma(X_t,S_t)$ depends on the current position of the particle, that is, on which boundary wheel it hits, and also on the current value of the spin $S_t$ at the time the boundary is hit. The spin of the particle is recorded by the process $S$. As we described it, it only updates when the particle is on the boundary, and we have chosen its amount of change to be linear with respect to the current spin, since this is the physically relevant situation. Indeed, angular momentum is conserved when two particles collide in absence of external interference. The spinning Brownian particle of our interest will collide against the revolving wheel and this system will locally maintain its total angular momentum. It is natural that the change of spin is given by a linear combination of the current spin  and the spin of the revolving boundary wheel [$\vec{g}(x)-\alpha(x)s$], also taking into account that part of the angular momentum is used in reflecting the particle in a non-normal direction (thus the factor $\alpha(x)$.) This model has inspired us to call the solution to \eqref{eq:sbm} spinning Brownian motion.

Even though our inspiration for the model comes from the spinning ball bouncing off of a moving boundary, from the mathematical point of view it is natural to regard the process $(X,S)$ as a multidimensional reflected difussion in $D\times\RR^p$ with degeneracy, due to the absence of a diffusive motion in the $p$ components of $S$. Setting $Z=(X,S)$ we can write \eqref{eq:sbm} as
\begin{align}
\label{eq:zsbm}
dZ_t = \sigma_0(Z_t)dB_t + \vec{\kappa}(Z_t)dL_t,
\end{align}
where 
$\sigma_0(x,s)$ is the $(n+p)\times n$ matrix  obtained from $\sigma(x)$ by augmenting it with zeroes, and $\vec\kappa(x,s) = \rpar{\vec\gamma(x,s),\vec g(x)-\alpha(x) s}$. Since we have $\partial(D\times\RR^p)=\partial D\times \RR^p$, the local times of \eqref{eq:sbm} and \eqref{eq:zsbm} are the same because $Z$ is in the boundary of its domain if and only if $X\in\partial D$, and the interior normal to $D\times\RR^p$ is just $(\vec n,0_p)$ where $0_p$ is the zero vector in $\RR^p$. 

Equation \eqref{eq:zsbm} does not fall within the domain of the submartingale problem of Stroock and Varadhan \cite{StV71} since the diffusion matrix $\sigma_0$ is not elliptic, and even though existence and uniqueness of a solution to equation \eqref{eq:zsbm} are direct to establish, their counterpart for the submartingale problem is  more subtle.

An alternative to the classical submartingale approach was introduced by Lions and Sznitman in  \cite{LiS84}, where existence of reflected diffusions driven by a general semimartingale was shown, but that result only holds for smooth, bounded domains. Their approach is based on an analytical solution to the deterministic Skorokhod problem (see also \cite{Cos92} and \cite{DuI93,DuI93c} for some non-smooth cases), but it does not yield many probabilistic results, such as the Feller property, that we need to study spinning Brownian motion and its stationary distribution.

In \cite{BaB10}, a reflected process with inert drift is studied. Existence is obtained by constructing a reflected Brownian motion, and then a drift is added through a Girsanov transformation. One could regard such process as a reflected diffusion with non-elliptic generator. Although their treatment of existence differs considerably from ours, the main ideas they use in the proof of uniqueness of the stationary distribution can be applied with some modification to our case.

Understanding the structure of the stationary distribution of spinning Brownian motion has been one of our main interests. To this end, in the second part of this article we restrict  to spinning Brownian motion processes with diffusive matrix $\sigma=I_n$, the identity matrix of size $n\times n$. First, we show that the spin $S_t$ eventually hits and stays within a certain compact, convex set, which is independent of the starting position of the process. A classical result for Feller processes then yields the existence of a stationary distribution. The most challenging part is to prove that spinning Brownian motion admits a unique stationary distribution under the following crucial assumption on the vector field $\vec{g}$:

\begin{description}
\item[A1] There are $p+1$ points $x_1,\ldots,x_{p+1}$ on the boundary of $D$ such that for every $y\in\RR^p$, there exist non negative coefficients $\lambda_j$ such that 
$y = \sum_{j=1}^{p+1} \lambda_j \vec g(x_j)$.
\end{description}

We start with an intermediate result that apparently has little to do with the stationary distribution, and it is interesting on its own. We identify the components of an exit system (see Section \ref{se:excursions}, or \cite{Mai75} for a definition) for excursions away from the boundary, in terms of the local time $L_t$ of the process, and a family of excursion measures $\h_x$ that has been constructed in the build up for Theorem 7.2 in \cite{Bur87}. It has been pointed out to us that it is possible to use the exit system $(L_t,\h_x)$ to construct a stationary distribution for the process $(X,S)$, in a similar manner as it is done in Theorem 8.1 in \cite{Get79}. We do not need to use such machinery to obtain a stationary distribution, as spinning Brownian motion happens to be a Feller process that stays within a fixed compact set, and thus existence of a stationary distribution follows from classic results.

Our proof of uniqueness of the stationary distribution is an adaptation, and somehow a generalization, of an analogous result for Brownian motion with inert drift, recently proved by Bass, Burdzy, Chen and Hairer in \cite{BaB10}. Although the literature on stationary distributions of elliptic reflecting diffusions is vast, we have not found other models or results for process that are similar to spinning Brownian motion. Nonetheless, some of the results available for elliptic reflected diffusions have been helpful to understand the challenges of our research. The reader can consult the articles \cite{HaL85,HaW87,HaW87exp}  to obtain an idea of the treatment of the problem in the case of an elliptic generator, and see the differences with our approach to the problem.

One result that helps  characterize the stationary distribution was developed by Weiss in his unpublished thesis \cite{Wei81}, by a test that involves only the candidate to stationary measure, the infinitesimal generator of the diffusion, and the vector field defining the directions of reflection. A recent extension of this result to some non-smooth domains was carried out by Kang and Ramanan in \cite{KaR14}. Both results ask for the submartingale problem associated to the diffusion to be well-posed, but they do not ask for ellipticity of the infinitesimal generator, and thus, they apply to our setting. We  make  use of this characterization of the stationary distribution to produce an explicit example of stationary distribution in a specific case of \eqref{eq:sbm}, and we also use this characterization to show that the stationary distribution of spinning Brownian motion never has a product form with factors related to $X$ and to $S$, under smooth conditions. We believe that this result is true in a more general setting, but our proof relies in the smoothness of the domain and of the coefficients.

\subsection{Notation}

We list some symbols that are used recurrently in the paper:  $\RR^q$ is the $q$-dimensional Euclidean space. For $A\sub\RR^q$, $\overline A$ denotes the closure of the set $A$, and $\partial A$ denotes its Euclidean boundary. The function $\ind_A$ is the indicator (characteristic) function of the set $A$. If  $x\in\RR^q$, its components are denoted $x^i$, with $i=1,\ldots,q$. Te canonical vector $\mathbf{e}_i\in\RR^q$ has all components equal to $0$, except for $\mathbf{e}_i^i=1$. For vectors $\vec a,\vec b\in\RR^q$ its inner product will be denoted by $\vec{a}\cdot\vec{b}$. Balls in $\RR^q$ are the sets $B(x,r)=\kpar{z\in\RR^q : \abs{z-x}<r}$, where the dimension $q$ will be clear from the context. The norm of $x\in\RR^q$ will be denoted by $\abs{x}$, and if $f$ is an $\RR^q$-valued function, we denote by $\norm{f}_\infty = \sup_x\abs{f(x)}$ its $L^\infty$ norm. 

The set of $q\times q$, positive definite matrices with real entries is denoted by $M_q^+(\RR)$.

When integrating, $dx$ refers to the Lebesgue measure differential, but the Lebesgue measure of a measurable set $A\sub\RR^q$ will be denoted by $m^q(A)$. Also, for $A,B\sub\RR^q$, the standard distance between these sets is denoted by $\mathrm{dist}(A,B)=\inf\kpar{\abs{x-y} : x\in A,y\in B}$. If $A=\kpar{x}$, we will write $\mathrm{dist}(x,B)$ for simplicity. Given a sufficiently smooth domain $G\subset\RR^q$, we denote the surface measure of its boundary $\partial G$ by $\nu(dx)$. This measure typically appears in the context of Green's theorem and harmonic measures.

The set of twice continuously differentiable functions in a domain $O$ is denoted by $C^2(O)$ or just $C^2$ when the domain $O$ is clear from the context. Similarly, We denote by $C^2_b(O)$ or $C^2_b$ its subspace of bounded functions. The space $C^{1,2}_{0}([a,b]\times O)$ correspond to the functions that have one continuous derivative in $(a,b)$, are continuous in $[a,b]$, and are twice continuously differentiable in $O$. The subscript $0$ indicates that the functions under consideration have compact support.

\newcommand{\D}{\mathcal{D}}

Let $(\O,\F,\PP)$ be a complete probability space and let $\kpar{B_t}_{t\ge 0}$ be a $q$-dimensional Brownian motion adapted to a filtration $\kpar{\F_t}_{t\ge 0}$, which satisfies the usual conditions. For a Borel set $G\sub \RR^q$, we denote its hitting time by $T_G=\inf\kpar{t>0 : B_t\in G}$. The Skorohod space will be denoted by $\D(q)=D([0,T];\RR^q)$, where balls are considered with respect to the Skorohod metric $d_\D$, and denoted by $B_\D(f_0,r)$. For further Probabilistic notation, we refer the reader to \cite{KaS91}.

\subsection{Outline}
\label{outline}

The paper is organized as follows: In Section \ref{se:main} we present the core results of our research: strong existence of \sbm, the excursion decomposition of paths, and many lemmas on the stationary distribution. We leave out two key (and harder) proofs for Section \ref{se:proofs}.

The submartingale problem characterization allows us to provide both explicit and numerical examples of the stationary distribution for \sbm\ processes. Section \ref{se:examples} contains applications of the main result of \cite{Wei81} and \cite{KaR14}.

\section{Main results}
\label{se:main}
\newcommand{\thetav}{\boldsymbol\theta}

\subsection{Existence of Spinning Brownian motion}
Our first results establish existence and uniqueness of spinning Brownian motion both as a solution to the SDE \eqref{eq:sbm}, and as a solution to a submartingale problem.

We recall our basic assumptions on the coefficients of equation \eqref{eq:sbm}: the domain $D\sub\RR^n$ is assumed to be of class $C^2$ and bounded. We assume that $\sigma$ is an $n\times n$ uniformly elliptic, and Lipschitz continuous. The field $\vec\gamma \colon \partial D\to\RR^n$: it is Lipschitz continuous, and for every $x\in\partial D$ we have $\vec\gamma(x)\cdot\vec n(x) =1$. The vector field $\vec{g}\colon\partial D\to\RR^p$ is Lipschitz, and $\alpha : \partial D\to\RR$ is a uniformly positive, Lipschitz continuous function. Notice that as $\overline D$ is compact, all continuous functions on $D$ are bounded.

\begin{theorem}
\label{th:existence_sde}
Under the conditions stated above, the stochastic differential equation with reflection
\begin{align*}
\left\{
\begin{array}{rl}
dX_t &= \sigma(X_t)dB_t + \vec{\gamma}(X_t,S_t)dL_t, \\
dS_t &= \spar{\vec{g}(X_t)-\alpha(X_t)S_t}dL_t,
\end{array}
\right.
\end{align*}
has a unique strong solution.
\end{theorem}
\proof The theorem follows almost immediately form Corollary 5.2 in \cite{DuI93} after a simple manipulation of the equation. Set $Z_t = (X_t,S_t)$,  and for $z=(x,s)\in\RR^n\times\RR^p$ define $\sigma_0(x,s)^T = \spar{\sigma(x)^T 0_{n,p}^T}$, where $0_{n,p}$ is the $n\times p$ matrix with all entries equal to zero. Also, set $\vec\kappa\rpar{x,s} = \begin{pmatrix} \vec\gamma(x,s) \\ \vec{g}(x)-\alpha(x)s\end{pmatrix}$. Equation \eqref{eq:sbm} can be rewritten as
\begin{align}
\label{eq:sbm_z}
dZ_t = \sigma_0(Z_t) dB_t + \vec{\kappa}(Z_t)dL_t,
\end{align}
where $Z_t\in \overline{D\times\RR^p}$. Note that $\partial (D\times\RR^p) = \partial D\times \RR^p$, and the interior normal to $D\times\RR^p$ is just the interior normal $\vec n$ to $D$ enlarged by $p$ zeros. It is straight forward to check that the local times in \eqref{eq:sbm} and \eqref{eq:sbm_z}
are equivalent, and so, these two equations are indeed equivalent.

Equation \eqref{eq:sbm_z} fits the framework of Corollary 5.2 in \cite{DuI93}, except for the fact that the domain of the reflected diffusion is unbounded, which makes the vector field $\vec\kappa$ unbounded. To fix this, we apply the Corollary 5.2 to equation \eqref{eq:sbm_z} in the domain $D\times B(0,2n)$ to obtain a process $Z^n_t$, and define for $m\in\NN$ the stopping times $\tau^n_m = \inf\kpar{t \geq 0 : \abs{S^n_t}>m}$. We apply Corollary 5.2 from \cite{DuI93}: for large $k\ge n\ge m$ we have that $Z^k_{t\wedge \tau^k_m}=Z^n_{t\wedge\tau^n_m}$ which implies that $\tau^k_m=\tau^n_m$ under $\PP_z$, whenever $\abs{z}<m$. In particular, this shows that  $\tau_m = \tau^{k}_m$ is well defined, and thus, under $\PP_z$, the process $Z_t=Z^k_{t\wedge\tau_m}$ is also well defined for $t<\tau_m$, for $\abs{z}<m$. 

It follows that a process $Z_t$ solving  \eqref{eq:sbm_z} can be defined up to time $\tau=\sup_n \tau_n$. We next show that $\tau=\infty$ under $\PP_z$ for any $z\in \overline D\times\RR^p$. Indeed, since the local time $L_t$ is continuous and of bounded variation, we have that the quadratic variation of $S^k$ is zero. Let $\alpha_0 = \inf\kpar{\alpha(x) : x\in\partial D}$. By It\^o's formula,
\begin{align*}
e^{\alpha_0 L^k_{t\wedge \tau_m}}\abs{S^k_{t\wedge \tau_m}}^2 &= \abs{S_0}^2 + \int_0^{t\wedge \tau_m} \!\!\! 2e^{\alpha_0 L^k_u}S^k_u dS^k_u + \alpha_0\int_0^{t \wedge \tau_m} \!\!\! e^{\alpha_0 L^k_u}\abs{S^k_u}^2 dL^k_u \\
&\leq \abs{S_0}^2 +  \int_0^{t\wedge \tau_m} e^{\alpha_0 L^k_u} \rpar{ 2\norm{\vec g}_\infty \abs{S^k_u} - \alpha_0\abs{S^k_u}^2} dL^k_u \\
&\leq \abs{S_0}^2 +  \int_0^{t\wedge \tau_m} e^{\alpha_0 L^k_u} \frac{\norm{\vec g}^2_\infty}{\alpha_0}dL^k_u \\
&= \abs{S_0}^2 +  \frac{\norm{\vec g}^2_\infty}{\alpha_0^2}\rpar{e^{\alpha_0 L^k_{t\wedge \tau_m}}-1}.
\end{align*}
Therefore,
$$
\abs{S^k_{t\wedge \tau_m}}^2 \leq \abs{S_0}^2e^{-\alpha_0L^k_{t\wedge \tau_m}} + \frac{\norm{\vec g}^2_\infty}{\alpha_0^2},
$$
which shows that for large enough $m$, we have $\tau_m=\infty$ under $\PP_z$. Otherwise, taking $t\to\infty$ above yields $m^2\leq \abs{S_0}^2 + \alpha_0^{-2}\norm{\vec g}^2_\infty$ for all large enough $m$, which is obviously a contradiction. This shows that the process $Z_t$ solves \eqref{eq:sbm_z} for all $t\ge 0$. 

Uniqueness follows from the same idea. Any two processes solving \eqref{eq:sbm_z} would coincide up to time $\tau_m$ by Corollary 5.2, \cite{DuI93}, thus they would coincide for all times.\qed

\begin{remark} The computation that lead to the bound in $\abs{S^k_{t\wedge\tau_m}}^2$ carries over to $S_t$. In this case we obtain $\abs{S_t}^2 \leq \abs{S_0}^2e^{-\alpha_0 L_t} + \alpha_0^{-2}\norm{\vec g}^2_\infty$. As it is deduced from the submartingale formulation we describe next, $L_t$ grows to infinity a.s., which implies that for large times the spin process lives in a neighbourhood of the ball $B(0,\alpha_0^{-1}\norm{\vec g}_\infty)$, and so any stationary distribution of $(X_t,S_t)$ must be supported at most in the closure of ${D\times B(0,\alpha_0^{-1}\norm{\vec g}_\infty)}$. For this reason, from this point on we will consider  $(X_t,S_t)$ as a bounded diffusion.
\end{remark}

One very successful way of constructing diffusion processes with boundary conditions was developed by Stroock and Varadhan in  \cite{StV71}. Their submartingale problem proved to be a successful extension of their ideas developed to treat the well-known martingale problem. The following survey on the submartingale problem is based on their original presentation. 

Let $G$ be a non-empty, open subset of $\RR^k$, such that: 
\begin{enumerate}[(i)]
\item there exists $\phi\in C^2_b(\RR^k;\RR)$ such that $G=\phi^{-1}(0,\infty)$, and $\partial G=\phi^{-1}(\kpar{0})$.
\item $\abs{\nabla\phi(x)}\geq 1$ for all $x\in\partial G$.
\end{enumerate}
The following functions will also be given:
\begin{enumerate}[(i)]
\item $a:[0,\infty)\times G\to M_n^+(\RR)$ which is bounded and continuous,
\item $b:[0,\infty)\times G\to \RR^k$ which is bounded and continuous,
\item $\vec\eta:[0,\infty)\times\partial G\to\RR^k$ which is bounded, continuous, and satisfies that ${\vec\eta(t,x)\cdot\nabla\phi(x)}\geq \beta >0$ for $t\geq 0$ and $x\in\partial G$.
\item $\rho : [0,\infty) \times \partial G \to [0,\infty)$ which is bounded and continuous.
\end{enumerate}
Define, for $u\geq 0$ and $x\in G$ 
\begin{align}
\label{eq:orbm_generator}
\L_u = \frac12 \sum_{i,j=1}^k a_{i,j}(u,x)\pd{^2}{x^i\partial x^j} + \sum_{i=1}^k b^i(u,x) \pd{}{x^i};
\end{align}
and, for $u\geq 0$ and $x\in\partial G$
$$
J_u = \sum_{i=1}^k \eta^i(u,x)\pd{}{x_i}.
$$
We say that a probability measure $\PP$ on $(\O,\F)$ solves the submartingale problem on $G$ for coefficients $a,b,\vec\eta$ and $\rho$ if $\PP\rpar{Z_t\in \overline{G}}=1$, for $t\geq 0$, and
$$
f(t,Z_t) - \int_0^t \ind_G(Z_u)\spar{\pd{f}{u} + \L_uf}(u,Z_u)\ du
$$
is a $\PP$-submartingale for any $f\in C_0^{1,2}\rpar{[0,\infty)\times\RR^k}$ satisfying 
$$
\rho \pd{f}{t} + J_t f \geq 0\qquad \mbox{on} \qquad [0,\infty)\times \partial G.
$$
We say the the submartingale problem is well-posed if it has a unique solution.

We next show that the SDE \eqref{eq:sbm} and its submartingale problem formulation are equivalent. This is done in order to access all the probabilistic results that the submartingale problem framework provides. In our case, the domain $D$ satisfies conditions (i) and (ii) above, which are easy to extend to $G=D\times\RR^p$. We set $a= \sigma_0  \sigma_0^T$, $b\equiv 0$ and $\vec\eta=\vec\kappa$, as in equation \eqref{eq:sbm_z}. In our case $\rho\equiv 0$.

\begin{theorem}
\label{eq:sde_submg}
The solution to \eqref{eq:sbm} constructed in Theorem \ref{th:existence_sde} is the unique solution to the associated submartingale problem.
\end{theorem}
\proof From It\^o's formula, it is direct to check that the solution to \eqref{eq:sbm_z} solves the submartingale problem. We only need to show uniqueness.

Let $Z^{*}_t$ be a solution to the submartingale problem. We are going to show that $Z^*_t$ is a weak solution to \eqref{eq:sbm_z}. 

Indeed, by Theorem 2.5 in \cite{StV71}
there exists an increasing, continuous process $t\mapsto L^*_t$ such that $dL^*_t$ is supported in the set $\kpar{Z^*_t\in\partial D\times\RR^p}$, and by Theorem 2.1 in \cite{StV71}, we have that for all $\thetav\in \RR^{n+p}$ the following is a martingale:
\begin{align}
\label{eq:ex_martingale}
M^{\thetav}_t=\exp\spar{ \thetav \cdot \rpar{Z^*_t - Z^*_0 - \int_0^t \vec\kappa(Z^*_u)dL^*_u}-\frac12\int_0^t \thetav \cdot a(Z^*_u)\thetav du }
\end{align}
For $i,j=1,\ldots,n$ let $M^{ij}$ be the martingale above, obtained from $\thetav=\eta\mathbf{e}_i+\lambda\mathbf{e}_j$, where $\eta,\lambda\in\RR$, and $\mathbf{e}_j$ is  the $j$-th vector in the canonical basis of $\RR^{n+p}$. Note that we only care about the $n$ first vectors in this basis. By doing a second-order Taylor expansion of $M^{ij}_t$ in the variables $\eta,\lambda$ we readily obtain that for $i,j=1,\ldots,n$,
\begin{align}
N^{j}_t =  {Z^*_t}^j-{Z^*_0}^j-\int_0^t \vec\gamma(Z^*_u)^j dL^*_u
\end{align}
is a continuous martingale with quadratic cross-variation given by
\begin{align}
\label{eq:cross01}
\ppar{N^i,N^{j}}_t =  \int_0^t [\sigma\sigma^T](Z^*_u)^{ij}du.
\end{align}
Since $\sigma$ is a bounded, elliptic matrix (here is crucial that $i,j\leq n$), we have that these cross-variation processes are absolutely continuous functions of $t$, thus, in view of Theorem 4.2 and Remark 4.3 in \cite{KaS91} we conclude that there is an $n$-dimensional Brownian motion $\kpar{W_t}$ in $\rpar{\Omega,\F,\F_t,\PP}$ and an $n\times n$ matrix valued, adapted process $\kpar{\chi_t}$, with
$$
{Z^*_t}^j-{Z^*_0}^j-\int_0^t \vec\gamma(Z^*_u)^j dL^*_u = \sum_{k=1}^n\int_0^t \chi^{j,k}_udW^k_u,\qquad j=1,\ldots,n.
$$
From \eqref{eq:cross01}, and It\^o's isometry it follows by continuity that for all $t>0$
$$
\chi_t \chi^T_t = \sigma(Z^*_t)\sigma(Z^*_t)^T.
$$
Since $\sigma$ is uniformly elliptic, the process $\Gamma_t=\sigma(Z^*_t)^{-1}\chi_t$ is well defined and so is $B_t = \int_0^t \Gamma_u dW_u$. It is easy to check that $\Gamma_t$ is unitary for all $t>0$, and that $\kpar{B_t}$ is a Brownian motion adapted to $\kpar{\F_t}$ by using Levy's Theorem. It follows that 
$$
N_t = \int_0^t \chi_u dW_u = \int_0^t \sigma(Z^*_u) \Gamma_u dW_u = \int_0^t \sigma(Z^*_u) dB_u,
$$
as desired.

By using $\thetav = \lambda\mathbf{e}_j$, with $j=n+1,\ldots,n+p$, and doing  a Taylor  expansion in the variable $\lambda$ in \eqref{eq:ex_martingale}, we readily see that for $Z^*_t = (X^*_t,S^*_t)$ we have that
$$
\rpar{{Z^*_t}^j - {Z^*_0}^j - \int_0^t \rpar{\vec{g}(S^*_u)-\alpha(X^*_u)S^*_u}dL^*_u}^m
$$
is a martingale starting from zero for every $m\in\NN$, and thus is identically zero. This completes the proof that $Z^*=(X^*,S^*)$ is a weak solution of \eqref{eq:sbm}, and thus the law of $Z^*$ must be the one of the unique solution to \eqref{eq:sbm}, completing the proof of the theorem.
\qed

The following representation formula simplifies the analysis of the spin process $S_t$.

\begin{lemma}
\label{le:spin}
Let $(X_t,S_t)$ solve equation \eqref{eq:sbm}. Define $Y_t = \exp\rpar{\int_0^t \alpha(X_u)dL_u}$. Then we have $dY_t = \alpha(X_t)Y_t dL_t$, and the spin process $S_t$ has the pathwise representation
\begin{align}
\label{eq:spin}
S_t = Y_t^{-1} S_0 + Y_t^{-1}\int_0^t \frac{\vec{g}(X_u)}{\alpha(X_u)} dY_u.
\end{align}
Also, the support of any stationary distribution of $(X_t,S_t)$ must be contained in the closure of $D\times H_{\vec{g},\alpha}$, where $H_{\vec{g},\alpha}$ is the convex hull of the set $\kpar{\frac{\vec{g}(x)}{\alpha(x)}: x\in\partial D}$.
\end{lemma}
\proof Since $t\mapsto L_t$ is increasing and continuous, $dL_t$ is a Riemann-Stieltjes measure and the first assertion is a consequence of the chain rule. The spin process has zero quadratic variation as $L_t$ does. We compute:
\begin{align*}
d(S_tY_t) &= \rpar{\vec{g}(X_t)-\alpha(X_t)S_t}Y_tdL_t + S_tY_t\alpha(X_t)dL_t = \vec{g}(X_t)Y_tdL_t,
\end{align*}
and so $d(S_tY_t) = \alpha^{-1}(X_t)\vec{g}(X_t)dY_t$. Since $Y_t\geq 1$ for all $t>0$, \eqref{eq:spin} follows from integration of the equation above, and division by $Y_t$.

\newcommand{\chull}{\overline H_{\vec{g},\alpha}}
Next we prove the assertion on the support of stationary distributions. Since $X_t\in\overline D$, it is enough to show that for any stationary distribution $\mu$ and open set $A\sub\RR^p\setminus \chull$, we have that $\mu(\overline D\times A)=0$. Moreover, it is enough to consider the sets  $A_n = \kpar{s\in\RR^p : \textrm{dist}(s,H_{\vec{g},\alpha})>\frac1n}$.

Let $C_t = (Y_t-1)^{-1}\int_0^t \frac{\vec{g}(X_u)}{\alpha(X_u)} dY_u$. It is clear that $C_t\in\chull $ for all $t>0$. It's not hard to arrive at the estimate
\begin{align*}
\textrm{dist}(S_t,H_{\vec{g},\alpha})\leq\abs{S_t -C_t} 
&\leq \rpar{\abs{S_0} +  \frac{\norm{\vec g}_\infty}{\alpha_0}}e^{-\alpha_0 L_t}.
\end{align*}
It follows that
\begin{align*}
\PP_{x,s} (S_t \in A_n) &\leq \PP_{x,s} \rpar{ L_t <
\alpha_0^{-1}\abs{\log n \rpar{\abs{s}+\frac{\norm{\vec g}}{\alpha_0}}}  }
\end{align*}
which converges to zero as $t\to\infty$, as $L_t\to\infty$ a.s. Since $\mu$ is a stationary distribution, we have for all $t>0$,
$$
\mu(\overline D\times A_n) = \int \PP_{x,s} (S_t \in A_n) \mu(dxds),
$$
and by dominated convergence, we deduce that $\mu(\overline D\times A_n)=0$ for all $n\in\NN$, as we wanted to show.\qed

\begin{remark}
\label{re:spinhull}
Equation \eqref{eq:spin} relates the process $S$ and the set $\chull$ in the following way: let $\P=\kpar{t_k}_{k=0}^N$ be a partition of the interval $[0,t]$. If the length of the longest interval of the partition is denoted $\abs{\P}$, then we have
$$
S_t = \lim_{\abs{\P}\to 0} S_0 Y_t^{-1} + Y_t^{-1} \sum_{k=0}^{N-1} \frac{\vec{g}(X_{t_k})}{\alpha(X_{t_k})} \rpar{Y_{t_{k+1}}-Y_{t_k}} = \lim_{\abs{\P}\to 0}S_0\lambda_N +  \sum_{k=0}^{N-1} \lambda_k \frac{\vec{g}(X_{t_k})}{\alpha(X_{t_k})},
$$
where $\lambda_N = Y_t^{-1}> 0$, and $\lambda_k = Y_t^{-1}\rpar{ Y_{t_{k+1}}-Y_{t_k} }\geq 0$. It is clear that $\sum_{k=0}^N \lambda_k = Y_t^{-1}\rpar{Y_t-Y_0+1} =1$, which means that $S_t$ is in the closure of the convex hull of $\kpar{S_0}\cup H_{\vec{g},\alpha}$. In particular, if $S_0\in \chull $ then $S_t \in  \chull$ for all $t\geq 0$.
\end{remark}

\subsection{Exit system for excursions away from the boundary}
\label{se:excursions}

From this section on, we restrict to the case $\sigma(x)=I_n$. We next introduce the notion of Exit System, first developed by Maisonneuve in \cite{Mai75}, although the following definitions are taken from \cite{BaB10}.

Let $Z$ be a standard Markov process taking values in a domain $E\sub\RR^k$ with boundary 
$\partial E$. We attach to $E$ a ``cemetery'' point $\Delta$ outside of $\overline E$, and we denote by $\C$ the set of   functions $g:[0,\infty)\to\RR^k\cup\kpar{\Delta}$ that are continuous in some interval $[0,\zeta)$ taking values in $\RR^k$, and are equal to $\Delta$ in $[\zeta,\infty)$.

A family of {excursion laws} $\kpar{\h_{z}}_{z\in\partial E}$, is a family of sigma-finite measures on $\C$ such that the canonical process is strong Markov on $(t_0,\infty)$ under $\h_z$, for every $t_0 > 0$, with the transition probabilities of the process $Z$ killed upon hitting $\partial E$. Moreover, $\h_z$ gives zero mass to paths which do not start from $z$. 

Excursions of $Z$ from $\partial E$ will be denoted $e$ or $e_s$, i.e, if $s<u$ and $Z_s,Z_u\in\partial E$, and $Z_t\notin \partial E$ for $t\in (s,u)$, then $e_s = \kpar{e_s(t) = Z_{t+s},\ t\in [0,u-s)}$ and the lifetime of such excursion is given by $\zeta(e_s)=u-s$. By convention, $e_s(t)= \Delta$ for $t \geq\zeta$. 

Let $L^*_t$ be an additive functional of $Z$, and let $\xi_t=\inf\kpar{s\geq 0 : L^*_s\geq t}$ the corresponding right inverse. Let $I$ be the set of left endpoints of all connected components of $(0,\infty)\setminus\kpar{t\geq 0 : Z_t\in\partial E}$. 

\begin{theorem}[Theorem 1 in \cite{Mai75}]\label{th:exit_formula} 
There exists a positive, continuous additive functional $L^*$ of $Z$ such that, for every $z\in\overline E$, any positive, bounded, predictable process $V$, and any universally measurable function $f:\C\to [0,\infty)$ that vanishes on  excursions $e_t$ identically equal to $\Delta$, 
\begin{align}
\label{exit}
\EE_{z}\spar{\sum_{t\in I} V_t f(e_t)} = \EE_z\spar{\int_0^\infty \!\!\! V_{\xi_s} \h_{Z(\xi_s)}(f) ds} = \EE_z\spar{ \int_0^\infty  \!\!\! V_t \h_{Z_t}(f) dL^*_t}.
\end{align}
Standard notation is used for $\h_z(f) = \int_{\C} f d\h_z$. 
\end{theorem}

The previous result is a specialised version of the exit system formula, and a pair $(L^*_t,\h_{z})$ satisfying \eqref{exit} is called an exit system.

The exit system formula provides a technical tool to reconstruct the process $Z_t$ excursion by excursion. A very nice use of the excursion formula allows us to ``count'' excursion with a given property. For instance, let $\Gamma$ be the set of excursions from $\partial E$ starting at a time $t\in [a,b)$ and going through an open set $U\sub E$. We set $V_t\equiv 1$ and $f=\ind_\Gamma$ in the exit formula to obtain
$$
\EE_z\spar{\sum_{a\le u < b} \ind_\Gamma(e_u)} = \EE_z\spar{\int_a^b \h_{Z_u}(\Gamma)dL^*_u}
$$
The left hand side is the expectation of the number of excursions in $\Gamma$, which can be computed using the exit system $(L^*_t,\h_z)$ according to the right hand side.

We remind the reader that a measurable set $B$ is called polar if $\PP_x(T_B<\infty)=0$ for all $x\in B^c$, otherwise the set $B$ is called nonpolar.  It is known that excursions laws can be picked in a standard way, precisely, there are unique (up to a multiplicative constant) excursions laws such that
\begin{enumerate}[(i)]
\item $\h_z\rpar{\lim_{t\to 0} Z_t \neq z} = 0$, for all $z\in\partial E$,
\item  $0 < \h_z(T_B <\infty) <\infty$ for all compact, nonpolar sets $B\sub E$. 
\end{enumerate}
For \sbm, note that the spin $S$ is constant during excursions, and since the dynamics of $X$ when inside $D$ do not depend of $S$, we see that excursion laws for $X$ should not depend on the value of $S$ at the beginning of the excursion.  Moreover, the law of $(X,S)$ killed upon hitting $\partial D\times\RR^p$ equals to that of Brownian motion killed upon hitting $\partial D$ times the law of $S_0$, so we can restrict our space of excursions to trajectories of the $X$-component of $(X,S)$ only. It follows that any excursion measure must be a product of the form $\h_{x,s}=\h_x\otimes\delta_{\kpar{s}}$, where $\h_{x}$ is an excursion measure representing paths of $X$ only. 

\begin{theorem}
\label{LHexit}
Let $\PP^D$ be the law of Brownian motion killed upon exiting $D$. For $x\in\partial D$, $s\in\RR^p$, define 
\begin{align}
\label{exitLaw}
\h_{x} \df \lim_{\e\downarrow 0} \frac{1}{\e}\PP^D_{x+\e\vec{n}(x)},\qquad \h_{x,s}\df\h_x\otimes\delta_{\kpar{s}},
\end{align} 
and let $L_t$ be the local time of $(X,S)$, satisfying equation \eqref{eq:sbm}, with $\sigma=I_n$. Then $\h_{x,s}$ is a sigma-finite measure, strongly Markovian with respect to the filtration of the driving Brownian motion $B_t$, and $(L_t,\h_{x,s})$ is an exit system from $\partial D\times\RR^p$ for the process $(X,S)$.
\end{theorem}
\proof See section \ref{se:proofs}.\qed

We will often disregard the dependence on $s$ of the measure $\h_{x,s}$, and will write $\h_x$ instead for the sake of simplifying the notation. 

Notice that the exit system formula does not offer a natural way to normalize the measures $\h_{x}$. Moreover, if $(A_t,\h_x)$ is an exit system and $\eta(x)$ defines a positive, measurable function in $E$, then $(\eta\cdot A_t,\eta(x)^{-1}\h_x)$ also defines an exit system. The excursion measures $\h_x$ introduced in the previous theorem have been used by Burdzy \cite{Bur87} to establish a canonical choice of an exit system for reflected Brownian motion in Lipschitz domains. 

For spinning Brownian motion, excursions from $\partial D$ start exactly at times when the local time increases, and thus it is natural that for some positive function $\eta:\partial D\to\RR$ we have that $(\eta\cdot dL_t,\h_x)$ is an exit system, because excursions of \sbm\ don't look any different from those of reflected Brownian motion. By reasoning as above, an exit system for \sbm\ should be $(dL_t,\eta^{-1}(x)\h_x)$. The theorem then proves that $\eta\equiv 1$.
\subsection{The stationary distribution}
\label{se:stationary}

The main goal of the section is to prove existence and uniqueness of the stationary distribution of spinning Brownian motion. Recall that we are in the case $\sigma(x)= I_n$.

One of the  issues with the diffusion $(X_t,S_t)$ is the lack of a driving Brownian motion for the coordinates related to the spin. At an intuitive level, this means that the spin $S_t$ could be confined to very small regions of the space, regions having Hausdorff dimension less than $p$, and consequently the support of the stationary distribution of the process could be singular with respect to Lebesgue measure. To make sure this is not the case, we  impose the following condition on the infinitesimal change of $S_t$, more precisely, on the function $\vec g$: 

\noindent \textbf{A1} There are $p+1$ points $x_1,\ldots,x_{p+1}$ on the boundary of $D$ such that for every $y\in\RR^p$, there exist non-negative coefficients $\lambda_j$ such that 
$y = \sum_{j=1}^{p+1} \lambda_j \vec g(x_j)$.

From now on, we assume that \textbf{A1} holds, and we fix the points $x_1,\ldots,x_{p+1}$ that realize it. Notice that if $y=\sum_{j=1}^{p+1}\vec g(x_j)$, then \textbf{A1} implies that $-y$ has an expansion with non-negative coefficients, and so we have that for every $\e>0$ there are coefficients $\eta_j>0$ such that $0 = \sum_{j=1}^{p+1} \eta_j\vec g(x_j)$, and $\sum_{j=1}^{p+1}\eta_j <\e$.

\begin{lemma}
\label{le:openNBHD}
The set $U_\e = \kpar{\sum_{j=1}^{p+1} \eta_j\vec g(x_j) : \eta_j\geq 0,\ 0<\sum_{j=1}^{p+1} \eta_j <\e}$ is an open neighborhood of zero for every $\e>0$. 
\end{lemma}
\proof
We claim that \textbf{A1} ensures that there are positive numbers $\eta^0_1,\ldots\eta^0_{p+1}$ such that $0=\sum_{j=1}^{p+1} \eta^0_j\vec{g}(x_j)$. Indeed, set $y_0=\sum_{j=1}^{p+1} \vec{g}(x_j)$, and use \textbf{A1} to pick non-negative $\eta^1_1,\ldots\eta^1_{p+1}$ such that $-y_0=\sum_{j=1}^{p+1} \eta^1_j \vec{g}(x_j)$. Then, set $\eta^0_j=1+\eta^1_j$.

Fix $\e>0$, and let $y = \sum_{j=1}^{p+1} \eta_j\vec g(x_j) $ for some $\eta_j\geq 0$, with $0< \sum_{j=1}^{p+1} \eta_j <\e$. Then, for every $\lambda>0$ we have that $y =  \sum_{j=1}^{p+1} \rpar{\eta_j + \lambda\eta^0_j}\vec g(x_j) $. The coefficients $\eta_j+\lambda\eta^0_j$ are all positive, and their sum is smaller than $\e$ for $\lambda$ small enough. It follows that
$$
U_\e =  \kpar{\sum_{j=1}^{p+1} \eta_j\vec g(x_j) : \eta_j > 0,\ 0<\sum_{j=1}^{p+1} \eta_j <\e},
$$
which is clearly open.\qed

\renewcommand{\vga}{\vec g/\alpha}
\renewcommand{\hull}{H_{\vec{g},\alpha} }

As we have already seen in Lemma \ref{le:spin}, the set of convex combinations of $\vga$  plays a significant role  in the characterization of the support of the stationary distribution of spinning Brownian motion. We named this set $\hull$, and we refer to it as the \emph{convex hull of $\kpar{\vga}=\kpar{\vec g(x)\alpha(x)^{-1} : x\in\partial D}$}.
In Remark \ref{re:spinhull}, we have seen that when started on $\hull$, the spin process $S_t$ lives forever in the closure of this set. 

We prepare to prove that the stationary distribution exists, is unique, and its support corresponds to the closure of $D\times\hull$. The proof consists of four steps. In the first one (Proposition \ref{pr:hitUe}), we use a support theorem and continuity results for the Skorohod map to show that for any given  point $z\in D$, $T>0$ and $\e>0$, the probability of $(X_T,S_T)$ to be in a ball of radius $\e$ around the final point $(z,0)$ is positive, no matter what the initial position is. In the second step, we use the results of Section \ref{se:excursions} and excursion theory to show how the path of $X_t$ can be decomposed into several excursions, and how spinning Brownian motion up to the first hitting time of a ball $U\sub D$ can be obtained from \sbm\ conditioned on never hitting $U$, adding a suitable ``last excursion'' that hits $U$. This construction is then used in the third step to patch together a spinning Brownian motion from several independent spinning Brownian motions $Y^j_t$. In the final step, we show how to condition each of the $Y^j$'s on hitting the 
boundary of $D$ only at certain places, and deduce  that the law of $S_t$ is bounded below by a measure which has a density with respect to Lebesgue measure. This procedure is detailed in the proof of  Theorem \ref{th:open0}.

\begin{lemma}
\label{SPbv}
Let $D$, $\vec{\tau}, \vec{g}$ and $\alpha$ be as above. Let $T>0$, and $z \in D$. Assume that \textbf{A1} holds. Then, for any $(x_0,s_0)\in\overline D\times\RR^p$ there is $\omega\in C([0,T];\RR^n)$ with bounded variation such that there is a unique $(x,s)\in C([0,T];\overline D\times\RR^p)$ satisfying $(x(0),s(0))=(x_0,s_0)$, $(x(T),s(T)) = (z,0)$, and for $t \in [0,T]$.
\begin{align*}
x(t) &= x_0 + \omega(t) + \int_0^t \vec{\gamma}(x(u),s(u)) dl(u), \\
s(t) &= s_0 + \int_0^t \vec{g}(x(u)) - \alpha(x(u))s(u)\ dl(u),
\end{align*}
Here, $l(\cdot)$ is a continuous and increasing function that only increases when $x(t)\in\partial D$, that is  $l(t) = \int_0^t \ind_{\partial D}(x(u))dl(u)$.
\end{lemma}

\proof We first show uniqueness. The idea is to apply Theorem 4.1 in \cite{LiS84} to the driving function 
$w_t= (x_0+\omega(t),s_0)$, and use the reflection field $\vec\kappa(x,s) = (\vec{\gamma}(x,s),\vec{g}(x)-\alpha(x)s)$.  Since by Remark \ref{re:spinhull}, the function $s(t)$ is bounded, there is no issue to apply Theorem 4.1 in \cite{LiS84} to our case.

Next we construct a function $\overline\omega\in  C([0,T];\RR^n)$ with bounded variation, and a solution $(x,s)$ of  the system above. Consider the uniform partition $0<a_1<b_1<a_2<\cdots<b_{p+1}<T$ of $[0,T]$.

To construct $\ow$ and the associated solution, set $\ow(0)=0$ and for $t\in (0,a_1]$, let $\ow(t)$ be defined as any continuous function with bounded variation such that $x_0+\ow(t)\in D$, and $x_0+\ow(a_1)=x_1$. It is clear that any solution $(x,s)$ has to satisfy $x(t)=x_0+\ow(t)$, $s(t)=s_0$, and $l(t)=0$ up to time $a_1$. Next, we want to keep $x(t)$ at $x_1$ from $a_1$ to $b_1$. In view of \eqref{eq:spin}, for $t\in [a_1,b_1]$ we set $y_1(t) = \exp \spar{\alpha(x_1)(l(t)-l(a_1))}=e^{\alpha(x_1)l(t)}$ and so
\begin{align*}
s(t) &= y_1(t)^{-1}s_0 + \vec g(x_1) y_1(t)^{-1} \int_{a_1}^{t} e^{\alpha(x_1)l(u)}dl(u) \\
&= y_1(t)^{-1} s_0 + \frac{\vec g(x_1)}{\alpha(x_1)}  y_1(t)^{-1} \spar{ y_1(t)-1}.
\end{align*}

By setting $l(t) = 0 + \eta_1(t-a_1)$ for $t\in [a_1,b_1]$, where $\eta_1\geq 0$ is to be determined, we obtain that both $l$ and $s$ are continuous for $t\le b_1$. All this implies that we need to define
$$
\ow(t) = x_1-x_0-\eta_1\int_{a_1}^t \vec{\gamma}\rpar{x_1, s(u)} du.
$$
Uniqueness in $[a_1,b_1]$ follows directly form the fact that the equation above defines a continuous function with bounded variation. Thus, the functions $(x,s)$ defined above correspond to the unique solutions to the Skorokhod problem for $\ow$ in $[a_1,b_1]$. 

Next, we want to keep $s(t)$ constant in $[b_1,a_2]$, while we move $x(t)$ from $x_1$ to $x_2$.  To this end, pick a curve with bounded variation $\zeta_1:[b_1,a_2]\to \overline D$ such that $\zeta_1(b_1)=x_1$, $\zeta_1(a_2)=x_2$ and $\zeta_1(t)\in D$ for other values of $t$. Set $l(t)=l(b_1)$ for $t\in [b_1,a_2]$,  and $\ow(t) = \zeta(t) - x_1 +\ow(b_1)$. It is clear that the only solution with bounded variation in this interval is $x(t) = \zeta_1(t)$ and $s(t) = s(b_1)$.

We iterate this process by keeping $x(t)$ at $x_j$ in $[a_j,b_j]$, and by defining $l(t) = l(b_{j-1})+\eta_j(t-a_j)$, $y_j(t) = \exp[\eta_j\alpha(x_j)(t-a_j)]$ in that interval. This way, the function $s(t)$ must satisfy
\begin{align}
\label{eq:k02}
s(t) &= [y_1(b_1)\cdots y_{j-1}(b_{j-1})y_j(t)]^{-1}s_0 \ + \\
\nonumber & \hspace{2cm}+ \sum_{m=1}^{j} \frac{\vec g(x_m)}{\alpha(x_m)} [y_m(b_m)\cdots y_{j-1}(b_{j-1}) y_j(t)]^{-1} \spar{y_m(b_m)-1}
\end{align}
for $t\in [a_j,b_j]$. The calculation leading to such equation, though tedious, is   straight-forward to carry out by splitting the integral in $[0,t]$ into integrals in the sets $[a_j,b_j]$ and $[b_j,a_{j+1}]$, and using our definition of $l(u)$ in each of those intervals. In the interval $[a_j,b_j]$, define $\ow(t)$ by
$$
\ow(t) = \ow(a_j) - \eta_j\int_{a_j}^t \vec\gamma (x_j,s(u)) du.
$$
Once again the unique solution in this interval for this $\ow$ is $(x_j,s(t))$.

For $t\in [b_j,a_{j+1}]$, we find a curve with bounded variation $\zeta_j$ going from $x_j$ to $x_{j+1}$ through $D$, and set $\ow(t) = \zeta_j(t) - x_j +\ow(b_j)$ and $l(t)=l(b_j)$. The unique solution is then $(\zeta_j(t),s(b_j))$. This procedure can be also done so that $x(T)=z$.

It remains to show that we can choose the values of $\eta_1,\ldots,\eta_{p+1}\geq 0$ such that $s(T)=0$. At time $T$ we find that
$$
s(T) = s_0\prod_{m=1}^{p+1}y_m(b_m)^{-1}  + \sum_{m=1}^{p+1} \frac{\vec g(x_m)}{\alpha(x_m)} [y_m(b_m)-1] \prod_{i=m}^{p+1} y_i(b_i)^{-1},
$$
and to obtain $s(T)=0$ we need 
$$
-s_0 = \sum_{m=1}^{p+1} \frac{\vec g(x_m)}{\alpha(x_m)} [y_m(b_m)-1] \prod_{i=1}^{m-1} y_i(b_i).
$$
By \textbf{A1} there are non-negative $\lambda_1,\ldots,\lambda_{p+1}$ such that $-s_0 = \sum_{m=1}^{p+1} \vec g(x_m)\lambda_m$, so we need to choose the numbers $\eta_m$ so that $y_m(b_m)$ satisfies $\alpha(x_m)\lambda_m = [y_m(b_m)-1]\prod_{i=1}^{m-1} y_i(b_i)$. This is easily achieved by an inductive procedure, and the lemma is proved.
\qed

\begin{proposition}
\label{pr:hitUe}
Let $D$, $\vec{\tau}$, $\vec{g}$, $\alpha$, $T>0$, and $z\in D$ as in Lemma \ref{SPbv}. For any $r>0$, there is $p_0>0$, such that for all $(x_0,s_0)\in \overline{D}\times {\chull}$, the inequality
$$
\PP_{x_0,s_0} \rpar{(X_T,S_T) \in B(z,r)\times B(0,r)} \geq p_0,
$$
holds.
\end{proposition}
\proof 
Let $\PP$ be the law of standard Brownian motion in $\RR^n$, and fix $(x_0,s_0)\in\overline{D}\times\RR^p$.  By pathwise uniqueness, we know that for a.e.  $\omega:[0,\infty)\to\RR^n$, there is a unique pair $(x,s)\in C([0,T];\overline D\times\RR^p)$, such that for $t\in [0,T]$
\begin{align*}
x(t) &= x_0 + \omega(t) + \int_0^t \vec{\gamma}(x(u),s(u)) dl(u), \\
s(t) &= s_0 + \int_0^t \spar{\vec{g}(x(u)) - \alpha(x(u))s(u)} dl(u),
\end{align*}
where, $l(\cdot)$ is a continuous and increasing, satisfying $l(t) = \int_0^t \ind_{\partial D}(x(u))dl(u)$, that is, it only increases when $x(t)\in\partial D$. It is standard to call this function $l(\cdot)$ the local time.

Let $\O$ be the set of continuous $\omega\in D([0,T];\RR^n)$ such that this uniqueness hold. We emphasize that $\PP(\O)=1$ and that the function $\ow$ constructed in Lemma \ref{SPbv} belongs to $\O$. Define $\Gamma$ in $\O$ by the assignment $x_0+\omega\mapsto(x,s)$ as above. We claim that $\Gamma$ is continuous at $\ow$, where continuity is taken in the sense of uniform convergence in compact sets. Indeed, let $\omega_j\in \O$ be a sequence converging uniformly in $[0,T]$ to $\ow$. Then, by setting $z(t) = (x,s)(t)$ and $\vec{\zeta}(z) = (\vec\gamma(x,s),\vec{g}(x))$, we have that the $\RR^{n+p}-$valued functions $(\omega_j,0)$ converge uniformly to $\overline\eta=(\ow,0)$. By Theorem 3.1 in \cite{Cos92}, we have that the unique solutions $(x_j,s_j)$ to the Skorokhod problem with reflecting vector $\vec\zeta$ in $\overline D\times\RR^p$, and corresponding driving function $(x_0+\omega_j,s_0)$ is relatively compact, and any limit is a solution of the corresponding problem with driving function $(x_0+\ow, s_0)$. By uniqueness ($\ow\in\O$), we deduce that $(x_j,s_j)\to (\overline x,\overline s)$ in $\D(n+p)=D([0,T],\RR^{n+p})$. But as all the involved functions are continuous, we actually deduce that the latter convergence is uniform in $[0,T]$.

In particular, there is $\delta>0$ 
such that if $(\omega,0)\in\O\cap B_{\D(n+p)}(\overline\eta,\delta)$, then the associated solution to the Skorohod problem $(x,s)\in B_\D((\overline x,\overline s),r)$. Thus we have
\begin{align*}
\PP_{x_0,s_0} \rpar{(X_T,S_T) \in B(z,r)\times B(0,r) } &\geq \PP_{x_0,s_0} \rpar{(x,s) \in B_{C[0,T]}((\overline x,\overline s),r) } \\
&\geq \PP_{x_0,s_0}\rpar{ (\omega,0)\in B_{\D(n+p)}(\overline\eta,\delta)} \\
&= \PP_{x_0}(\o \in B_{\D(n)}(\overline{w},\delta)),
\end{align*}
The last term is greater than some positive constant $p_0$ that depends on $x_0, s_0$, and $r>0$, by the support theorem of Brownian motion. The Feller property of $\PP$, and a standard compactness argument applied in $\overline{D}\times\chull$, let us choose $p_0$ independently of $x_0, s_0$.
\qed

\begin{corollary}
\label{co:tauFinite}
Let $r>0$ and $\tau=\inf\kpar{t>0 : S_t\in B(0,r)}$. Then $\tau $ is finite almost surely.
\end{corollary}
\proof 
Let $N(b)$ be the event ``$S_t$ is not in $B(0,r)$ for any $t\in [0,b]$''. Then $\PP_{x,s}\rpar{\tau<\infty} = 1 - \lim_{b\to\infty}\PP_{x,s}\rpar{N(b)}$, where the limit is clearly  decreasing. 

By Proposition \ref{pr:hitUe} we have that $\PP_{x,s}(N(T))\leq 1-q$, and a standard application of the Markov property shows that
%
$\PP_{x,s}\rpar{N(nT)}\leq (1-q)^n$, which yields $\PP_{x,s}\rpar{\tau<\infty}=1$.
\qed

We next proceed to introduce some results about the stationary distribution of spinning Brownian motion. Our method is very much an adaptation of the proof of Theorem 6.1 in \cite{BaB10}. Such argument involves a decomposition of the law of $X_t$ into several reflecting processes that are somewhat independent of each other. To the reader familiar with excursion theory, ``independence'' is achieved by using suitable exit systems. This decomposition allows us to control both the local time and the trajectory of the process before hitting a fixed open set $U$, and deduce that no stationary measure can be null in $U$.

\begin{theorem}
\label{th:open0}
Let $Z=(X,S)$ be spinning Brownian motion solving \eqref{eq:sbm} with $\sigma\equiv I_n$, and fix $z\in D$. There is $r>0$, $t_0>0$, a constant $c_{10}>0$, and an open set $V_*\sub \hull$, such that for every initial condition $(x,s)\in B(z,r)\times B(0,r)$, the law of the random variable $Z_{t_0} =(X_{t_0},S_{t_0})$ is bounded below by a measure with a density with respect to $(n+p)$-dimensional Lebesgue measure. Moreover, this density is bounded below by $c_{10}$, on $D\times V_*$.  
\end{theorem}
\proof See section \ref{se:uniqueness_proofs}.

\begin{corollary}
\label{co:uniqueSD}
Spinning Brownian motion has a unique stationary distribution, supported in the closure of $D\times\hull$.
\end{corollary}
\proof 
Fix $T>0$, $z\in D$ as in Lemma \ref{SPbv}. Fix $r,t_0,c_{10}$ and $V_*$ as in Theorem \ref{th:open0}. Using $r>0$, pick $p_0>0$ as given by Proposition \ref{pr:hitUe}.

From Remark \ref{re:spinhull}, we know that any stationary distribution of \sbm\ has to be supported in the closure of $D\times \hull$, a bounded set, therefore we deduce that \sbm\ has at least one stationary distribution from the standard theory of Feller processes (see Theorem IV.9.3 in \cite{KuE86}). Let $\mu$ be one of them. For any open set $O\sub D\times V_*$,
\begin{align*}
\mu(O) &= 
\int _{\overline{D}\times\chull} \PP_{x,s}\rpar{(X_{T+t_0},S_{T+t_0})\in O} \ \mu(dxds) \\
&= \int _{\overline{D}\times\chull} \EE_{x,s}\rpar{ \PP_{X_{T},S_{T}} \rpar{(X_{t_0},S_{t_0})\in O}} \ \mu(dxds) \\
&\geq \int _{\overline{D}\times\chull} \EE_{x,s}\rpar{\ind_{B(z,r)\times B(0,r)}(X_T,S_T) \PP_{X_{T},S_{T}} \rpar{(X_{t_0},S_{t_0})\in O}} \ \mu(dxds) \\
&\geq c_{10}m^{n+p}(O)  \int _{\overline{D}\times\chull} \PP_{x,s}\rpar{(X_T,S_T) \in B(z,r)\times B(0,r) }\ \mu(dxds) \\
&\geq c_{10}m^{n+p}(O) p_0\mu\rpar{\overline{D}\times\chull} \\
&= p_0 c_{10} m^{n+p}(O),
\end{align*}
which means that any stationary distribution contains $D\times V_*$ in its support. If there were more than one stationary distribution, Birkhoff's ergodic theorem \cite{Sin94} implies that at least two of them, say $\mu$ and $\nu$, must be singular with respect to one another. But this contradicts the fact that both measures have $D\times V_*$ in their supports.
\qed


\section{Proofs}
\label{se:proofs}
\subsection{Proof of Theorem \ref{LHexit}}
\label{se:excursions_proofs}

%
The process $X_t$ behaves as Brownian motion inside $D$, so the exit laws of $X$ and the Brownian motion $B$ are the same. This suggests that we can use the same standard excursion measures for $X$ as for reflected Brownian motion, which are studied in \cite{Bur87}.
The fact that $\h_x$ is sigma-finite for all $x\in\partial D$, and strongly Markovian is proved in Theorem 7.2 in \cite{Bur87}. 

Let $(L^*_t, \h_x)$ be an exit system for $(X,S)$, where $L^*$ is the additive functional from Theorem \ref{th:exit_formula}. We will prove that it is possible to replace $L^*_t$ by the local time $L_t$ from equation \eqref{eq:sbm}. 

Fix $T>0$, $x_0\in\partial D$, and small enough $\e>0$ such that both $ B(x_0,\e)\cap D$ and $\Gamma=\overline{B(x_0,\e)}\cap\partial D$ are connected sets. The set $\partial\Gamma = \kpar{x\in\partial D : \abs{x-x_0}=\e }$ has surface measure zero, and since the surface measure and the harmonic measure are mutually absolutely continuous, almost surely no excursions of $X$ have ending points in $\partial\Gamma$. In particular, $\PP_{(x,s)}(T_{\partial\Gamma} < \infty)=0$ for all $x\in D$.

Let $G$ be a subdomain of $D$ such that $\overline{G}$ is a compact subset of $D$, and $\overline G\cap \overline{B(x_0,\e)} = \emptyset$. As usual $T_E$ denotes the hitting time of the set $E$ by the process $X$. Define $\A_G = \kpar{T_G < T_{\partial D}}$. 

Define the following sequences of stopping times: $\eta_0 = T_\Gamma$, and inductively set $\tau_{k+1} = \eta_k+T_{D\setminus B(x_0,\e)}\circ\theta_{\eta_k}$, and $\eta_{k+1} = \tau_{k+1}+ T_\Gamma\circ\theta_{\tau_{k+1}}$. We claim that $\eta_k$ grows to infinity almost surely. Assume that this is not the case, and choose $\lambda>0$ such that $\PP\rpar{\sup_{ k\geq 0} \eta_k \leq \lambda} = c_\lambda >0$. On the event $\kpar{\sup_{ k\geq 0} \eta_k \leq \lambda}$, the sequences $\eta_k$ and $\tau_k$ converge to the same random time $\tau\leq\lambda$. Since $X_{\eta_k}\in {\Gamma}$, and $X_{\tau_k}\in \overline{D}\setminus B(x_0,\e)$, it follows by continuity of paths that $X_\tau\in\partial\Gamma$. Thus, $T_{\partial\Gamma}<\infty$ on this event, which is a contradiction.

We will focus on excursions that start from $\Gamma$ before time $T>0$, and belong to $\A_G$. For $z=(x,s)$, $x\in \overline D$, we define
\begin{align*}
I_z(T,\Gamma;G) &= \EE_z\rpar{\sum_{t<T}\ind_{\Gamma}(X_t)\ind_{\A_G}(e_t)}.
\end{align*}
Let $L^*$ be an additive functional such that $(\h_x,L^*_t)$ is an exit system for excursions starting from $\partial D$, with $\h_x$ as in \eqref{exitLaw}. Using the exit system formula and the fact that $\eta_k,\tau_k\to \infty$, we get
\begin{align*}
I_z(T,\Gamma;G) &= \EE_z\rpar{ \int_0^T \ind_\Gamma(X_u) H_{X_u}(\A_G) dL^*_u} \\
&= \EE_z\rpar{ \sum_{k=0}^\infty \ind_{[0,T]}(\eta_k) \int_{\eta_k}^{\tau_{k+1}}\ind_{\Gamma}(X_u) H_{X_u}(\A_G) dL^*_u}.
\end{align*}
For $\eta_k\leq u <\tau_{k+1}$, it is clear that $X_u\in\partial D$ implies that $X_u\in \Gamma$. Thus, the term $\ind_{\Gamma}(X_u)$  is redundant on the right hand side of the last equality. By the strong Markov property
\begin{align}
\label{eq:proof_00} 
I_z(T,\Gamma;G) 
&= \EE_z\rpar{ \sum_{k=0}^\infty \ind_{[0,T]}(\eta_k) \EE_{(X_{\eta_k},S_{\eta_k})} \spar{\int_{0}^{\tau_{1}} H_{X_u}(\A_G) dL^*_u}} \\
\label{eq:proof_01}
&= \EE_z\rpar{ \sum_{k=0}^\infty \ind_{[0,T]}(\eta_k) \EE_{(X_{\eta_k},S_{\eta_k})} \spar{ \sum_{t<\tau_1} \ind_{\A_G}(e_t) }},
\end{align}
where in the last equality we have used the exit time formula with $V_t = \ind_{[0,\tau_1]}(t)$ and $f=\ind_{\A_G}$.

Notice that starting from $(y,v)\in\Gamma\times \RR^p$, the random variable $\sum_{t<\tau_1} \ind_{\A_G}(e_t)$ counts the amount of excursions in $\A_G$ (starting from $\Gamma$) that start before the process $X$ exits the ball $B(x_0,\e)$. This can happen only once, namely, for the excursion $e_t$ with $t=\sup_{u<\tau_1} \kpar{X_u\in\Gamma}$. Moreover, it happens exactly once on the event $\A_G\circ\theta_{\tau_1}$, so we have that 
\begin{align}
\label{eq:proof_02}
\EE_{(y,v)}\rpar{\sum_{t<\tau_1} \ind_{\A_G}(e_t)} = \PP_{(y,v)}\rpar{\A_G\circ\theta_{\tau_1}} = \EE_{(y,v)}\rpar{ \PP_{X_{\tau_1}}(\A_G)},
\end{align}
by the strong Markov property. Define $h(x) = \PP_x(\A_G)=\PP_x^D(\A_G)$ for $x\in \overline{D}\setminus {G}$. Since $\A_G$ depends only on the behaviour of $X_t$ up to $T_{\partial D}$, and $X_t$ is a Brownian motion inside $D$, a standard argument shows that $h$ is harmonic in $D\setminus\overline{G}$, with $h\vert_{\partial D}=0$ and $h\vert_{\partial G}=1$. Also, notice that $\nabla h(x)\cdot\vec{\gamma}(x,s) = \nabla h(x)\cdot\hat{n}(x)$ for $x\in\partial D$, since $h$ is constant along the boundary of $D$. Further, the definition of directional derivative shows that $ \nabla h(x)\cdot\hat{n}(x)=\h_x(\A_G)$. It follows that
$$
h(X_{t\wedge T_G}) - \int_0^{t\wedge T_G}  \h_{X_u}(\A_G) dL_u
$$
is a martingale under $\PP_{(y,v)}$, for $y\in \overline{D}\setminus G$. If $y\in\Gamma$, we have $\PP_{(y,v)}(\tau_1<T_G)=1$, and since $h$ is bounded, we have by the optional sampling theorem
$$
\EE_{(y,v)}\rpar{ \PP_{X_{\tau_1}}(\A_G)} = \EE_{(y,v)}\rpar{h(X_{\tau_1})} = \EE_{(y,v)}\rpar{ \int_0^{\tau_1} \h_{X_u}(\A_G)  dL_u }.
$$
Using this in \eqref{eq:proof_02}  and  \eqref{eq:proof_01} we obtain,
\begin{align*}
I_z(T,\Gamma;G) &= \EE_z\rpar{ \sum_{k=0}^\infty \ind_{[0,T]}(\eta_k) \EE_{(X_{\eta_k},S_{\eta_k})} \spar{\int_{0}^{\tau_{1}} H_{X_u}(\A_G) dL_u}}.
\end{align*}
Notice that last equation is just \eqref{eq:proof_00} with $L_t$ instead of $L^*_t$. Unravelling the steps that led to \eqref{eq:proof_00} gives
\begin{align}
\EE_z\rpar{ \int_0^T \ind_\Gamma(X_t) H_{X_u}(\A_G) dL^*_u} &= \EE_z\rpar{ \int_0^T \ind_\Gamma(X_t) H_{X_u}(\A_G) dL_u}.
\end{align}

A standard argument involving the monotone class theorem, shows that this last equation is not only valid for  $f=\ind_\Gamma$, but also for all bounded, measurable functions $f:\partial D\to\RR$. Standard estimates for Brownian motion show that the function $x\mapsto \h_{x}(\A_G)$ is bounded away from zero (see \cite{Bur87}),  so we can take $f(x) = e^{-\alpha t} \h_x(\A_G)^{-1}$ to get
\begin{align}\label{eq:exit_02}
\EE_{z} {\int_v^t e^{-\alpha t} dL^*_u} = \EE_{z} {\int_v^t e^{-\alpha t}  dL_u},
\end{align}
for any  positive $\alpha$, and arbitrary  $v<t$.   We can extend \eqref{eq:exit_02} in the following way: for $N\in\NN$ and $T>0$, define simple functions by
$$
f_{N,T}(t) = \ind_{\kpar{0}}(t) + \sum_{k=0}^N \ind_{\left(\frac{kT}{N},\frac{(k+1)T}{N}\right]}(t) \ e^{-\alpha \frac{(k+1)T}{N}}.
$$ 
It follows from \eqref{eq:exit_02} that
$
\EE_{z} {\int_0^\infty f_{N,T}(u) dL^*_u} = \EE_{z} {\int_0^\infty f_{N,T}(u)  dL_u}.
$
It is clear that $f_{N,T}(u)$ increases to $e^{-\alpha u}\ind_{[0,T]}(u)$ for all $T>0$, and so, by the monotone convergence theorem, we obtain  
$$
\EE_{z} {\int_0^\infty e^{-\alpha u} dL^*_u} = \EE_{z} {\int_0^\infty e^{-\alpha u} dL_u},
$$
that is, $L^*$ and $L$ have the same $\alpha$-potential functions. Since both $L^*$ and $L$ are continuous, it follows that $L^*=L$ a.s. by Theorem 2.13, Chapter 4, \cite{BlG68}. This shows that $(L_t,\h_x)$ is an exit system, and the theorem is proved.
\qed

\subsection{Proof of Theorem \ref{th:open0}}
\label{se:uniqueness_proofs}

%
We closely follow part of the proof of Theorem 6.1 in \cite{BaB10}. To simplify the notation, call $Z=(X,S)$, and we use the standard nomenclature $T^X_U$ for first hitting time of a set $U$ by the process $X$, and $\sigma_t = \inf\kpar{s\geq 0 : L_s\geq t}$ for the right inverse of local time.

We proceed to describe an exit system for a different, though related, process $X'$. Let $x_0\in D$ and $r>0$ be arbitrary but fixed, so that $\overline {B(x_0,r)} \sub D$. Set  $U= B(x_0,r)$ and let $Z'=Z^{U}$ be the process $Z$ conditioned by the event $\kpar{T^X_U > \sigma_1}$. It follows from Proposition \ref{pr:hitUe} and the strong Markov property that for any starting point in $D$, the probability of $\kpar{T^X_U > \sigma_1}$ is greater than zero. It is easy to see that $(Z'_t,L_t)$ is a time homogeneous Markov process. To be consistent with the notation, we will write $(Z'_t,L'_t)$ instead of $(Z'_t,L_t)$.
%
%
We will construct this exit system on the basis of $(L_t,\h_{x})$ because of the way that $Z'$ has been defined in relation to $Z$. It is clear that $L'$ does not change within any excursion interval of $X'$ away from $\partial D$, so we will assume that $\h'_{x,l}$ is a measure on paths representing $X'$ only. 
For $l\ge 1$ we let $\h'_{x,l} =\h_x$. Let $\widehat\PP^D_y$ denote the distribution of Brownian motion starting from $y\in D\setminus \overline U$, conditioned to hit $\partial D$ before hitting $U$, and killed upon exiting $D$. For $l<1$, we have 
\begin{align}
\label{hxl}
\h'_{x,l}(A) =  \lim_{\e\downarrow 0} \frac{1}{\e} \widehat\PP^D_{x+\e\hat n(x)} (A).
\end{align} 

Let $A_\star\sub \C$ be the event that the path of $X'$ hits $U$. It follows from the definition of $\h_x$ in Theorem \ref{LHexit} and \eqref{hxl} that for $l<1$,
\begin{align}
\label{hxl01}
\h'_{x,l}(A) = \h_x(A\setminus A_\star).
\end{align} 
One can deduce easily from standard estimates for Brownian motion that for some constants $c_3,c_4>0$ and all $x\in\partial D$,
\begin{align}
\label{hstar}
c_3 < \h_x(A_\star)< c_4.
\end{align}
Let $\sigma'_t = \inf \kpar{s \geq 0 : L'_s\geq t}$. The exit system formula \eqref{exit} and \eqref{hxl01} imply that we can construct $Z$ using $Z'$ as a building block, in the following way. Suppose that $Z'$ is given. We enlarge the probability space, if necessary, and construct a Poisson point process $\E$ with state space $[0,\infty)\times \C$ (see the definition of excursions) whose intensity measure conditional on the whole trajectory $\kpar{X'_t,\ t\geq 0}$ is given by
\begin{align}
\label{poisson}
\mu\rpar{[a,b]\times F} = \int_{1\wedge a}^{1\wedge b} \h_{X'_{\sigma'_t}}(F\cap A_\star) dt.
\end{align}
Since $\mu([0,\infty]\times \C) > c_3 $, the Poisson process $\E$ may be empty; that is, if the Poisson process is viewed as a random measure, then the support of that  measure may be empty. Consider the case when it is not empty and let $K_1$ be the minimum of the first coordinates of points in $\E$. Note that there can be only one point $(K_1,\e_{K_1})\in\E$ with first coordinate $K_1$, because of \eqref{hstar}. By definition of $\E$, we have that $\e_{K_1}=e_{\sigma'_{K_1}}$, where $e$ denotes excursions of $X'$. By convention, let $K_1=\infty$ if $\E=\emptyset$. Recall that $T^X_U = \inf \kpar{t>0 : X_t\in U}$ and let 
\begin{align*}
T_U^{X'} &= \inf \kpar{t>0 : X'_t\in U}, \\
T_\star &= \sigma'_{K_1} + \inf \kpar{t>0 : e_{K_1}(t)\in U}.
\end{align*}
It follows from the exit system formula that the distribution of the process
$$
\widehat{Z}_t = \begin{cases}
Z'_t & \text{if } 0\leq t \leq T^{X'}_U \wedge \sigma'_{K_1}, \\
\rpar{\e_{K_1}(t-\sigma'_{K_1}),S_{ \sigma'_{K_1}}} 
& \text{if } \E\neq\emptyset \text{ and } \sigma'_{K_1} < t \leq T_\star,
\end{cases}
$$
is the same as the distribution of $\kpar{Z_t, 0\leq t\leq T^X_U}$. So that we can refer easily to this construction later, we define $\Gamma(Z,U)=\widehat{Z}$.

Let $U_j = B(a_j,r)$ for $j=1,\ldots,p+1$, where $a_1=z$ was fixed in the statement, and $a_j\in D$ and $r>0$ are chosen so that the closure of the sets $U_j$ are pairwise disjoint and their union is a subset of $D$.

Let $Y^1=(X^1,S^1)$ be the process $\Gamma(Z,{U_2})$ satisfying that $X^1_0=x\in U_1$ and $S^1_0=s\in B(0,r)$. We will later shrink the value of $r$ if necessary. The process $Y^1$ is a \sbm\ starting from $(x,s)$, and observed until the first hit of $U_2\times\RR^p$, at time $T_1=\inf\kpar{t>0 : X^1_t\in U_2}$. Pick a \sbm\ $Z^2$, independent of $Y^1$, such that $Z^2_0= (X^{(2)},S^{(2)})$, where $X^{(2)}$ is uniform in $U_2$ and $S^{(2)}\sim S^1_{T_1}$, and define $Y^2=(X^2,S^2)$ as the process $\Gamma(Z^2,U_3)$, with $T_2=\inf\kpar{t>0 : X^2_t\in U_3}$. We have that $Y^2$ is  a \sbm\ observed up to the first time its $X$-component hits $U_3$.

Inductively, once $Y^1,\ldots,Y^j$ are defined, pick a \sbm\ $Z^{j+1}$ independent of $\kpar{Y^1,\ldots,Y^j}$, with initial distribution $Z^{j+1}_0=(X^{(j+1)},S^{(j+1)})$ where $X^{(j+1)}$ is uniform in $U_{j+1}$ and $S^{(j+1)}\sim S^j_{T_j}$, and define $Y^{j+1}=(X^{j+1},S^{j+1})$ as the process $\Gamma(Z^{j+1},U_{j+2})$, observed up to time $T_{j+1}=\inf\kpar{t>0:X^{j+1}_t\in U_{j+2}}$.
 It should be clear that the processes $Y^j$ are  independent, and have the distribution of a \sbm\  whose spatial component $X$ starts with uniform distribution in $U_j$, and the process is observed up to the first time it hits the set $U_{j+1}\times\RR^p$.

We will next use these processes $Y^j$ to patch together a \sbm\  $Z^*$. Note that for some $c_5 > 0$ and all $x,y \in U_{j+1},\  j = 1,\ldots,p+1$, 
$$
\PP_{x,\cdot} \rpar{X_1\in dy\ \text{and}\ X_t \notin \partial D\ \text{for}\ t\in [0,1]} \geq c_5 dy.
$$
We can assume that all $Z^j$'s and $Y^j$'s are defined on the same probability space. The last formula and standard coupling techniques show that on an enlarged probability space, there exist spinning Brownian motions $W^j$, $j = 1,\ldots,p+1$, with the following properties. For $1 \le j \le p$, $W^j_0 = Y^j_{T_j}$, and for some $c_6>0$
\begin{align}
\label{csix}
\PP\rpar{ W^j_1 = Y_0^{j+1}\ \text{and}\ W^j_t\notin \partial D\times\RR^p\ \forall\ t\in[0,1] \Big\vert \{Y^k\}_{k=1}^j,\ \{W^k\}_{k=1}^{j-1}} \ge c_6 	.		
\end{align}


The process $W^j$ does not depend otherwise on $\kpar{Y^k}_{k=1}^{p+1}$ and $\kpar{W^k}_{k\neq j}$. We define $W^{p+1}$ as a spinning Brownian motion  with $W^{p+1}_0 = Y^{p+1}_{T_{p+1}}$ but otherwise independent of $\kpar{Y^k}_{k=1}^{p+1}$ and $\kpar{W^k}_{k=1}^{p}$.

Let 
$$
F_j = \kpar{W^j_1 = Y^{j+1}_0 \ \text{and}\ W^j_t\notin \partial D \ \text{for}\ t\in [0,1]}.
$$
We define a process $Z^*$ as follows. We let $Z^*_t = Y^1_t$ for $0 \le t \le T_1$. If $F_1^c$ holds, then we let $Z^*_t = W^1_{t-T_1}$ for $t\geq T_1$. If $F_1$ holds, then we let $Z^*_t = W^1_{t-T_1-1}$ for $t\in [T_1,T_1+1]$ and $Z^*_t = Y^2_{ t-T_1-1}$ for $t\in [T_1+1,T_1+1+T_2]$. We proceed by induction. Suppose that $Z^*_t$ has been defined so far only for $t\in [0,T_1+1+T_2+1+\cdots+T_k]$, for some $k<n$. If $F_k^c$ holds then we let 
$$
Z^*_t = W^k_{t-T_1-1-T_2-1-\cdot-T_k}
$$
for $t\geq T_1+1+T_2+1+\cdots+T_k$. If $F_k$ holds, then we let 
$$
Z^*_t = W^k_{t-T_1-1-T_2-1-\cdot-T_k}
$$
for $t\in [T_1+1+T_2+1+\cdots+T_k,T_1+1+T_2+1+\cdots+T_k+1]$, and
$$
Z^*_t = Y^{k+1}_{t-T_1-1-T_2-1-\cdots-T_k-1}
$$
for $t\in [T_1+1+T_2+1+\cdots+T_k+1,T_1+1+T_2+1+\cdots+T_k+1+T_{k+1}]$. We let
$$
Z^*_t = W^{p+1}_{t-T_1-1-T_2-1-\cdots-T_{p+1}}
$$
for $t\ge T_1+1+T_2+1+\cdots+T_{p+1}$.

\bigskip

This construction makes the process $Z^*$ a spinning Brownian motion starting from $(x,s)$. To conclude this proof, we show that with a positive probability, the process $S$ can have ``almost'' independent and ``almost'' linearly independent increments over disjoint intervals of time. This is used to show that a conditioned version of $S$ has a density, or equivalently, that the law of $S$ is bounded below by a mesure with a density, in an appropriate open set.

Let $x_1, \ldots, x_{p+1}$ be points in $\partial D$ satisfying assumption \textbf{A1}. Since the matrix $[\vec g(x_1)|\cdots|\vec g(x_{p+1})]$ has rank $p$, it is possible to eliminate a column from it and still have a matrix with rank $p$. It follows, without loss of generality, that the vectors $\vec g(x_1),\ldots,\vec g(x_p)$ are linearly independent. For $1\leq j\leq p$, consider the set $C_j= \kpar{w\in\RR^n : \text{the angle between the vectors} \ \vec g(x_j) \text{ and } w \text{, is at most } \delta_0}$, for some $\delta_0>0$ so small that for any $z_j\in C_j$, $j=1,\ldots,p$, the vectors $\kpar{z_j}$ are still linearly independent. Let $\delta_1 > 0$ be so small that for every $j = 1,\ldots, p$, and any $x \in\partial D\cap B(x_j,\delta_1)$, we have $\vec g(x) \in C_j$.

Let $L^j$ the local time of $X^j$ on $\partial D$ and $\sigma^j_t = \inf \kpar{s\ge 0 : L_s^j \ge t}$. It is not hard to see that for some $p_2>0$, the probability that for every $j=1,\ldots,p$ we have $X^j_t \notin \partial D \setminus B(x_j,\delta_1)$, for $t \in [0,\sigma^j_t]$, is greater than $p_2$. Let 
$$
R_j = \sup \kpar{ t < T_j : Y^j_t \in \partial D } \quad \text{and} \quad Q_j =L^j_{R_j}.
$$

Consider the event $F_\star$ containing all trajectories such that for $j=1,\ldots, p$, we have $X^j_t\notin \partial D\setminus B{}(x_j,\delta_1)$ for $t \in [0,\sigma^j_1]$ and $R_j < \sigma^j_1$. The construction of the process $\Gamma(Z,U)$  shows  that for some constant $c_7>0$ the inequality $\PP_{z_0}(F_\star)\geq p_2 (1-e^{-c_7})^{p}$ holds.

Let $E^j_t = \exp \int_0^t \alpha(X^j_s) dL^j_s\geq 1$, and define the following collection of random variables
\begin{align*}
S^1(t_1,\ldots,t_p) &=  \rpar{ E^1_{\sigma^1_{t_1}} \cdots E^p_{\sigma^p_{t_p}} }^{-1} \rpar{s+  \int_0^{t_1} \alpha\rpar{X^1_{\sigma^1_u}}^{-1}\vec g\rpar{X^1_{\sigma^1_u}} E^1_{\sigma^1_u} du}, \\
S^j(t_j,\ldots,t_p) &=  \rpar{ E^j_{\sigma^j_{t_j}} \cdots E^p_{\sigma^p_{t_p}} }^{-1}  \int_0^{t_j} \alpha\rpar{X^j_{\sigma^j_u}}^{-1}\vec g\rpar{X^j_{\sigma^j_u}} E^j_{\sigma^j_u} du,\qquad j\geq 2.
\end{align*}
Notice that if $F_\star$ holds, then $S^j(t_j,\ldots,t_p)\in C_j$ for all  $t_j\in (0,1]$ and
$j=2,\ldots, p$. For $j=1$, one should notice that linear independency is stable under perturbations, thus, if $\abs{s}$ is small enough, we also have that $S^1(t_1,\ldots, t_p)\in C_1$ for $t_1\in (1/4,1]$, and $t_j\in (0,1]$ for all $j=2,\ldots, p$ . By shrinking the value of $r>0$ if necessary, we assume this holds. 

For any $0 \le a_k < b_k \le Q_k$ for $k=1,\ldots,p$, define the random set 
$$
\Lambda \rpar{ [a_1,b_1],\ldots,[a_{p},b_{p}] } = \kpar{ \sum_{j=1}^p S^j(t_j,\ldots,t_d): t_{k} \in [a_k,b_k]\ }.
$$ 
It is not difficult to estimate the $p$-dimensional volume of $\Lambda([a_1,b_1],\ldots,[a_p,b_p])$ by using the definition of $C_j$'s. First, by continuity, it follows that under $F_\star$ there is a positive constant $q$ such that  $(1-q\delta_0)e^{\alpha(x_j) t_j}\leq E^j_{\sigma^j_{t_j}}\leq (1+q\delta_0) e^{\alpha(x_j)t_j}$ for all $1\leq j\leq p$. By definition of the set $C_j$, it follows that for some  positive constant $\beta$,
$$
(1-\beta\delta_0) \frac{\vec g(x_j)}{\alpha(x_j)} \spar{e^{\alpha(x_j) t_j}-1} \leq \int_0^{t_j} \alpha(X^j_{\sigma^j_u}) ^{-1}\vec{g}(X^j_{\sigma^j_u}) E^j_{\sigma^j_u} du \leq (1+\beta\delta_0) \frac{\vec g(x_j)}{\alpha(x_j)} \spar{e^{\alpha(x_j)t_j}-1},
$$
where the inequality holds component by component. Define a function  $\vec v:[1/4,1]\times [0,1]^{p-1}\to\RR$ by
$$
\vec v(t_1,\ldots,t_p) = \sum_{j=1}^p {e^{-\sum\limits_{k=j}^p \alpha(x_k)t_k}} \frac{\vec g(x_j)}{\alpha(x_j)}\spar{e^{\alpha(x_j)t_j}-1}.
$$
It follows from the inequalities in this  paragraph that for some  constant $\eta >0$, independent of $\delta_0$, we have
$$
(1-\eta\delta_0) \vec{v}(t_1,\ldots,t_p) \leq \sum_{j=1}^p S^j(t_j,\ldots, t_p) \leq(1+\eta\delta_0) \vec{v}(t_1,\ldots,t_p),
$$
where the inequalities hold by components. If $\delta_0$ is small enough (and so is $r$) so that the vectors in different $C_j$'s are always linearly independent, the inequality above implies that there exist a constant $c_8$ independent of $a_k,b_k$ such that 
$$
c_8^{-1} \leq \frac{m^p\rpar{\Lambda([a_1,b_1],\ldots,[a_p,b_p])} }{ m^p\kpar{\vec v(t_1,\ldots,t_p) : t_k\in [a_k,b_k]} } \leq c_8.
$$
To compute the volume of the set in the denominator, in Lemma  \ref{le:jacobian} we calculate the Jacobian of $\vec v(\cdot)$, obtaining  $\det {D\vec v} =C \exp\rpar{-\sum_{k=1}^p k \alpha(x_k)t_k}$, which, as all variables $t_k$ are bounded,  readily yields that the $p$-dimensional volume of the random set $\Lambda([a_1,b_1],\ldots,[a_p,b_p])$ is bounded above by $c_4(b_1-a_1)\cdots(b_p-a_p)$ and below by $c_5 (b_1-a_1)\cdots(b_p-a_p)$.

Let us consider the processes $Z^*$ defined previously, conditioned on the sigma field 
$$
\G = \sigma \rpar{ \kpar{S^j(t_j,\ldots,t_p) ,\ t_1\in[1/4,1],\ t_2,\ldots,t_p \in [0,1] }_{j=2}^{p} }.
$$
The construction of $\Gamma(Z,U)$ makes the random variable $Q_j$ (defined earlier) the time component of a time-excursion Poisson random variable with variable intensity given by \eqref{poisson}. By \eqref{hstar}, it follows that $Q_j$ has a density with respect to the Lebesgue measure in $[0,1]$ that is bounded below.
In view of our remarks on the volume of $\Lambda$, it follows that conditional on $\G$, there is $c_9>0$ such that for all $s\in B(0,r)$, the vector 
$$
S(Q_1,\ldots,Q_p) = S^1(Q_1,\ldots,Q_p) + S^2(Q_2,\ldots,Q_p) +  \cdots + S^{p}(Q_{p})
$$ 
has a density with respect to the $p$-dimensional Lebesgue measure that is bounded below by $c_9$ on an open ball $U_*\subseteq \kpar{\vec{v}(t_1,\ldots,t_p) : t_j\in (1/4,1)}$. By Remark  \ref{re:spinhull}, we can assume that $U_*\sub \hull$. It is important to remark here that $c_9$ depends on $r$, but not on $s\in B(0,r)$. The same is true for $U_*$.

We can now remove the conditioning on $F_\star$ and conclude that the law of $S(Q_1,\ldots,Q_p)$ is bounded below by a measure with a density with respect to $p$-dimensional Lebesgue measure, and this density is bounded below on $U_*$. Define $E^*_t = \exp(\int_0^t \alpha(X^*_u) dL^*_u)$ and $S^*_t = {E^*_t}^{-1}s + {E^*_t}^{-1}\int_0^t \vec{g}(X^*_s)E^*_s d L^*_s$ and $T_\star = \sum_{j=1}^p T_j$, where $ L^*$ is the boundary local time of the spinning Brownian motion $X^*$. Our construction, and Lemma \ref{le:spin} show that $(X^*,S^*)$ is a \sbm\ starting from  $(x,s)$. Using conditioning on $F_\star$, we see that the law of $S^*_{T_\star}$ is bounded below on $U_*$ by the measure $c_9 m^p(\cdot \cap U^*)$.

The previous argument can be modified to show that for some fixed $t_0>0$, the law of the random variable $S^*_{t_0}$ is bounded below by a constant times the $p$-dimensional Lebesgue measure on a non-empty, open ball $V_*$. All  we need to do is, for small $\e>0$, find times $t_j>0$ such that  $T_j\in (t_j-\e,t_j+\e)$, with uniformly (in $j$) positive probability $q_{\e}$, and then further condition the stitched  process $X^*$ to satisfy $T_j\in (t_j-\e,t_j+\e)$. Set $t_*=\sum_{j=1}^p t_j$. This way, $T_* = \sum_{j=1}^p T_j \in (t_*-\e p,t_*+\e p)$ and since $X^*_{T_*}\in U_p$, and $U_p$ is away from the boundary, we can condition on $X^*$ to not to hit $\partial D$ in $[t_*-\e p,t^*+\e_p]$ and thus have $S^*_{t_*} = S^*_{T_*}$. We then choose  $t_0=t_*$. Since $X^*$ behaves as Brownian motion in $D$, we deduce that there is a constant $c_{10}>0$, such that the law of $(X^*_{t_0},S^*_{t_0})$ is bounded below on $D\times V_*$ by the measure $c_{10}m^{n+p}\rpar{\cdot\cap D\times V_*}$.
This holds for all initial conditions $(x,s)\in B(z,r)\times B(0,r)$.
\qed

\begin{lemma}
\label{le:jacobian} Let $\vec v(\cdot)$ be the function defined in the last step of the proof of Theorem \ref{th:open0}. There exists a constant $C\neq 0$, depending only on the vectors $\vec g(x_k)$ such that
$$
\det D\vec v(t_1,\ldots,t_p) = C \exp{\sum\limits_{k=1}^p- k \alpha(x_k)t_k} 
$$
for all $t_1,\ldots,t_p\in (0,1]$.
\end{lemma}
\proof 
A stright forward calculation shows that for $i=1,\ldots, p$:
\begin{align*}
\pd{\vec v}{t_i} 
%
%
&= e^{-\sum\limits_{k=i}^p  \alpha(x_k)t_k} \vec g(x_i) -\alpha(x_i)\sum_{j=1}^{i-1} e^{-\sum\limits_{k=j}^p  \alpha(x_k)t_k} \frac{\vec g(x_j)}{\alpha(x_j)}  \spar{ e^{\alpha(x_j)t_j} -1} .
\end{align*}

Let $\lambda_{j,i}$ the coefficient of $\vec{g}(x_j)$ in the expansion of $\pd{\vec v}{t_i}$ above. Because the vectors $\vec g(x_k)$ are linearly independent, these numbers are well defined. We highlight the fact that $\lambda_{j,i}=0$ for $j>i$. Let $T_{\vec{g}}$ the $p\times p$ matrix whose $j-$th column is $\vec g(x_j)$ and let $\Lambda$ be the matrix whose $(j,i)$ component is $\lambda_{j,i}$. The calculation above then simply says that $D\vec v = T_{\vec{g}} \Lambda$. Therefore, as $\Lambda$ is triangular
$$
\det D\vec v(t_1,\ldots,t_p) = \det T_{\vec g} \cdot \prod_{i=1}^p \lambda_{i,i} = \det T_{\vec g} \cdot \exp \rpar{ - \sum_{i=1}^p \sum_{k=i}^p \alpha(x_k)t_k }.
$$ 
The double sum on the right hand side equals  $\sum_{k=1}^p k \alpha(x_k)t_k$, by Fubini's theorem, and the proof is complete.
\qed

\section{Examples}
\label{se:examples}

In his unpublished thesis \cite{Wei81}, Weiss develops a test to characterize any invariant measure of the solution to a well-posed submartingale problem. His test only works in a smooth setting, and has been recently extended to a larger class of non-smooth domains by Kang and Ramanan in \cite{KaR14}. Even though we just need Weiss' result given our assumptions on the domain $D$, the fact that the result is available for more general ones opens a research line that we had not considered before. We lay out these results next.

Set $\vec\kappa(x,s)=\rpar{\vec\gamma(x,s),\vec{g}(x)-\alpha(x)s}$. Since $S_t$ is bounded in the stationary regime we can regard the vector $\vec\kappa$ as bounded. It follows that the following theorem from the unpublished dissertation of Weiss \cite{Wei81} can be applied to our setting: Let $\L = \frac12 \sum_{i,j=1}^n a_{i,j}(x)\pd{^2}{x_i\partial x_j} + \sum_{i=1}^n b_i(x) \pd{}{x_i}$ be a second order differential operator, where $a_{i,j}$ and $b_i$ are bounded, Lipschitz functions. Assume that a bounded, Lipschitz vector field 
$\vec\kappa$ is given on the boundary of a $C^2(\RR^d)$ domain $G$, such that $\vec\kappa\cdot\hat n(x)\ge\beta >0$ for $x\in\partial G$. Let $\phi$ be a $C^2(\RR^d)$ function defining the boundary of $G$.

\begin{theorem}[from \cite{Wei81}]
\label{th:weiss}
Let $\overline G$ be  compact in $\RR^d$ and $b_j$ and $\vec\kappa$ as before, suppose $\rpar{a_{i,j}(x)}$ is bounded, continuous, and positive semidefinite satisfying $\nabla\phi(x)^Ta(x)\nabla\phi(x)>0$  for $x$ in a neighborhood of $\partial G$ (i.e. the diffusion has nonzero random component normal to the boundary). Suppose that the submartingale problem for $a,b$ and $\vec\kappa$ is uniquely solvable starting from any $x\in \overline G$. 

Then, a probability measure $\mu$ on $\overline G$ is invariant for the diffusion if and only if $\mu(\partial G)=0$ and 
\begin{align}
\label{eq:weiss}
\int_G \L f(x) \mu(dx) \leq 0
\end{align}
for all $f\in C_b^2(\overline G)$ with $\nabla f\cdot\vec\kappa(x)\ge 0$  for $x\in \partial G$. 
\end{theorem}

This theorem has been successfully used by Harrison, Landau and Shepp \cite{HaL85} to give
an explicit formula for the stationary distribution of obliquely reflected Brownian motion in planar domains, in two cases: (a) the domain is of class $C^2(\CC)$ and bounded, and the reflection coefficient $\kappa$ has a global extension to a $C^2_b(\RR^2)$ vector field; and (b) the domain
is a convex polygon, and the reflection coefficient is constant in each face. Their technique to obtain an explicit representation is to assume that the stationary distribution has a density $\rho$ with respect to Lebesgue measure in the domain, integrate \eqref{eq:weiss} by parts to obtain a PDE with boundary conditions for $\rho$, and solve such equation. Our approach to obtain the stationary distribution for some specific cases of spinning Brownian motion is based on the same idea.

\subsection{Spinning Brownian motion in a wristband}
Consider spinning Brownian motion in the strip $\RR\times [-1,1]$ with
coefficients given by 
%
$g(y) = \frac12(\kappa+\beta)y + \frac12(\kappa-\beta)$ for positive constants $\kappa,\beta$, and the and $\tau(x,y;s) = \lambda\hat x\ind_{\kpar{1}}(y)$, where $\hat{x}=\mathbf{e}_1$ and $\lambda$ is a constant. We set $\alpha(x,y)=1$ for simplicity. The associated spinning Brownian motion $(X^x_t,X^y_t,S_t)$
solves the equation
$$
\left\{
\begin{array}{ll}
dX^y_t &= dB^y_t + \rpar{\ind_{\kpar{-1}}(X^y_t) - \ind_{\kpar{1}}(X^y_t)}dL_t ,\\ 
dX^x_t &= dB^x_t +\lambda\ind_{\kpar{1}}(X^y_t)dL_t, \\
dS_t &= \spar{g(X^y_t)-S_t} dL_t.
\end{array}\right.
$$

Note that the normal depends only on the $y$-coordinate, and so $X^y_t$ has the
distribution of reflected Brownian motion in $[-1,1]$. In particular, $L_t$
depends exclusively on $B^y_t$. Also, if we identify the points $x$ and
$x+2\pi$, the domain becomes a compact space and the existence of a unique
stationary distribution follows with minor and obvious modifications from our theorem. We will use this identification in what follows.

It is clear from the
equations that the law of $(X,S)$ starting from $(x,y,s)$ is the same as the law
of $(x+{X^*}^x,{X^*}^y,S^*)$, where $(X^*,S)$ starts from $(0,y,s)$. It follows that the stationary distribution is invariant under translations in the
$x$-coordinate. Thus, the stationary distribution of $(X^x,X^y,S)$ can be obtained from that of $(X^y,S)$ multiplying by $(2\pi)^{-1}dx$.

\begin{proposition}
The stationary distribution for the process $(X^x,X^y,S)$ is given by the density function $\rho(x,y;s)=a(s)y+b(s)$, where
\begin{align}
\label{eq:a_b}
a(s) = \displaystyle\frac{2}{\kappa+\beta}  \frac{s - \frac12 (\kappa-\beta) }{
\sqrt{(\kappa-s)(\beta+s)} }  \qquad\qquad b(s) = \displaystyle \frac{1}{\sqrt{(\kappa-s)(\beta+s)} }.
\end{align}
\end{proposition}
\proof It is enough to show that $\tilde\rho(y,s)=a(s)y+b(s)$ is stationary for the process $(X^y,S)$ which is a diffusion that solves a well posed submartingale problem in the domain $G=  (-1,1) \times (-\beta,\kappa)$. Set up this way, Theorem \ref{th:weiss} does not apply as $G$ is not of class $C^2$. Nonethless, since the spin process is bounded to $\overline{H}_{g}=[-\beta,\kappa]$ it is not hard to find  a bounded domain $ G_1$ of class $C^2(\RR^2)$ such that
$$
G\sub  G_1 \sub (-1,1)\times\RR,
$$
and apply the theorem in $ G_1$ for the density $\tilde\rho' = \tilde\rho\ind_{G}$. Another option is to use the recent version of Theorem \ref{th:weiss} for non-smooth domains due to Kang and Ramanan \cite{KaR14}.

Set $\vec\kappa(y,s)=(\vec n(y),{g}(y)-s)$, where $\vec n =\ind_{\kpar{-1}} - \ind_{\kpar{1}}$ is the negative of the sign function. The process $(X^y,S)$ that uniquely solves the submartingale process associated to  $\partial_{y,y}$ with boundary condition $\nabla f \cdot \vec\kappa(y,s)\geq 0$ for $y\in \partial {G_1}$ and $f\in C^2_b(\overline{{G}}_1)$. Since $\tilde\rho' =0$ outside of $G$, the following computation is straight forward by integration by parts: 
\begin{align*}
\int_{G_1} \partial_{yy} f(y,s) \tilde\rho'(y,s) dyds &= \int_{-\beta}^\kappa\int_{-1}^1  \partial_{yy} f(y,s) \tilde\rho(y,s) dyds 
= \int_{-\beta}^\kappa {\partial_{y} f} (y,s) \tilde\rho(y,s) - f(y,s)a(s)\Big\vert_{-1}^1 ds 
\end{align*}
The boundary condition $\nabla f\cdot\vec\kappa(y,s)\geq 0$ for $y=\pm 1$ translates into
\begin{align*}
 \spar{g(y)-s}\partial_s f (y,s) \geq \mathrm{sgn}(y)\partial_y f(y,s) & \qquad\qquad y=\pm 1.
\end{align*}
As $\tilde\rho\ge 0$, we see that $\displaystyle \partial_y f(y,s)\tilde\rho(y,s)\vert_{-1}^1 \leq \mathrm{sgn}(y)\spar{g(y)-s}\tilde\rho(y,s)\partial_s f(y,s)\vert_{-1}^1$. Also, direct computation shows that $\spar{g(y)-s}\tilde\rho(y,s) = \frac{2\mathrm{sgn}(y)}{\kappa+\beta}\sqrt{(\kappa-s)(\beta+s)}$  at $y=\pm 1$. Notice that this term vanishes both at $s=-\beta$ and $s=\kappa$, and its partial derivative with respect to $s$ equals to $-\mathrm{sgn}(y)a(s)$ at $y=\pm 1$. Doing integration by parts in $s$, and using the facts above we obtain:
$$
\int_{G'} \partial_{yy} f(y,s) \tilde\rho'(y,s) dyds \le 0,
$$
as desired.\qed

\

%

\subsection{No product form for the stationary distribution}

Let's assume that the domain $D$, and the coefficients $\vec\gamma, \vec{g},$ and $\alpha$ are smooth. We treat $\overline D\times \RR^p$ as an $(n+p)$-manifold with boundary $\partial D\times\RR^p$. All the vectors fields defined in the boundary can be extended in a smooth way to $D\times \RR^p$, for instance, fixing $s\in\RR^p$, and extending its components harmonically to $D$. Since the setting is smooth, standard results show that this extension is smooth as well.

For fixed $s\in\RR^p$, it is possible to find the flow $\theta^s(t,x)$ of the extension of $\vec{n}$. Working in local coordinates, it is easy to see that for a point $x_0\in\partial D$ the integral line $t\mapsto\theta^s(t,x_0)$ only intersects the boundary at $t=0$. A standard compactness argument then shows that there is an $\e>0$ such that the map 
$$
(t,x_0,s)\in [0,\e)\times\partial D\times \chull\mapsto(\theta^s(t,x_0),s)
$$ 
is a parametrization of a neigborhood $U$ of $\partial D\times  \chull$.

Let $h$ be a smooth function in $U$, and let $f_0$ be a smooth function in $\partial D\times\RR^p$. With the parametrization defined above, we have that
$$
f(\theta^s(t,x_0),s) = \int_0^t h(\theta^s(u,x_0),s) du + f_0(x_0,s)
$$
defines a smooth function $f$ in $U$ satisfying that
\begin{align}
\label{eq:bdryextension}
\nabla_x f(x,s) \cdot\vec{n}(x) = h(x,s),\qquad f(x,s)=f_0(x,s)\qquad x\in\partial D, s\in \overline H_{\vec{g},\alpha}.
\end{align}
Notice that if $\vec{\tau}(x,s)$ is a vector field tangential to $\partial D$, the product $\nabla_x f(x,s)\cdot \vec{\tau}(x,s)$ only depends on the boundary values of $f(x,s)$, that is, on $f_0(x,s)$. Thus, $\nabla_x f(x,s)\cdot \vec{\tau}(x,s) = \nabla_x f_0(x,s)\cdot \vec{\tau}(x,s)$.

\begin{proposition}
\label{pr:no_product}
Let $(X_t,S_t)$ be spinning Brownian motion with diffusion matrix $I_n$, and let the domain $D$ and the coefficients $\vec\gamma,\ \vec{g}$ and $\alpha$ be smooth. Then, the stationary distribution of $(X,S)$ does not have a product form.
\end{proposition}
\begin{proof}
Assume that the stationary distribution is $\mu(dx)v(ds)$. By Theorem \ref{th:weiss} we have that $\mu(\partial D)=0$. Since $X$ is a Brownian motion inside $D$,  the marginal $\mu$ admits a harmonic density $u(x)dx=\mu(dx)$.

Set $c_u=\int_{\partial D} u(x)\alpha(x) \nu(dx)$, $\vec{g}_u=c_u^{-1}\int_{\partial D} \vec{g}(x)u(x) \nu(dx)$, and $f_0(s)=\vec{g}_u\cdot s-\frac12\abs{s}^2$. Notice that $f_0$ does not depend on $x$. Also, set $h(x,s) = -\nabla_s f_0(s)\spar{\vec{g}(x)-\alpha(x)s}$, and construct the function $f$ satisfying \eqref{eq:bdryextension}, as we did before stating the proposition. It is clear that $f\in C^2(\overline D\times \chull)$ and that 
\begin{align*}
\nabla f(x,s)\cdot\vec{\kappa}(x,s) &=\nabla_xf(x,s)\cdot\vec{n}(x)+\nabla_xf(x,s)\cdot\vec{\tau}(x,s)+\nabla_s f(x,s)\cdot [\vec{g}(x)-\alpha(x)s] \\
&= h(x,s) + \nabla_xf_0(s)\cdot\vec{\tau}(x,s)+\nabla_s f_0(s)\cdot [\vec{g}(x)-\alpha(x)s] =0.
\end{align*}

By Theorem \ref{th:weiss} we must have $\int_{D\times \hull} \Delta_x f(x,s) \mu(dx)v(ds)\leq 0$. Since $u$ is harmonic in $D$, and $f$ does not depend on $x$ on $\partial D$, by Green's formula and the definitions above, we have
\begin{align*}
\int_{D} \Delta_x f(x,s) u(x)dx &= -\int_{\partial D} u(x)\nabla_x f \cdot\vec{n}(x) \ \nu(dx) = 
\nabla_sf_0(s) \int_{\partial D} {u(x)\vec{g}(x)-s\alpha(x)u(x)} \ \nu(dx) \\
&= \nabla_sf_0(s) [c_u\vec{g}_u-c_us] = c_u\abs{\vec{g}_u-s}^2.
\end{align*}
It follows that 
\begin{align*}
\int_{D\times \hull} \Delta_x f(x,s) u(x)dx v(ds) &= 
c_u\int_{\hull} \abs{\vec{g}_u-s}^2 v(ds)
\end{align*}
is non-positive if and only if $\nu$ is a delta at $\vec{g}_u$, which contradicts the fact that $\mu$ is bounded below by a measure with density  in an open set of $D\times \hull$ (see Theorem \ref{th:open0}). This finishes the proof.
\end{proof}

\section*{Acknowledgements}
Most of the research leading to this publication was conducted during my graduate studies at the University of Washington. I would like to thank my advisor Krzysztof Burdzy for introducing me to the problem  discussed in this article, and for many conversations that helped me to successfully conduct this research. During this time, my research was funded by the NSF grant number DMS 090-6743. 

The preparation of this manuscript was partially funded by FONDECYT, project n$^\circ$ 3130724. We also acknowledge support of Programa Iniciativa Cientifica Milenio grant number NC130062
through the Nucleus Millenium Stochastic Models of Complex and Disordered Systems.

The author thanks an anonymous referee  for their detailed revision and helpful comments  that helped improved the exposition of this work, and motivated the author to improve and correct one of the proofs.
%

\section*{References}


\def\cprime{$'$}

\end{document}